\documentclass[titlepage,11pt]{article}
\usepackage{latexsym}
\usepackage{amsfonts}
\usepackage{amsmath}
\usepackage{float}
\usepackage{tikz}
\usetikzlibrary{decorations.pathmorphing}
\usetikzlibrary{positioning,shapes,arrows,decorations.markings}
\usetikzlibrary{decorations.pathreplacing,calligraphy}
\usetikzlibrary{positioning,arrows.meta}
\usetikzlibrary{fit}
\usepackage{etoolbox}

\newcommand*{\halfway}{0.5*\pgfdecoratedpathlength+.5*3pt}

\def\Stealtharrow{{\arrow[xshift=2pt+3.2\pgflinewidth]{Stealth[scale=1.3]}}}

\tikzset{
    midarrow/.style={decoration={markings,mark=at position #1 with {\Stealtharrow}},postaction={decorate}},
    midarrow/.default=0.5
}

\tikzset{->-/.style 2 args={decoration={
  markings,
  mark=at position \halfway with {\arrow[rotate=#2]{Stealth[scale=1.3]}}},postaction={decorate}},
  ->-/.default = {0.55}{0}
}
\tikzset{
dot/.style = {circle, fill, minimum size=#1,
              inner sep=0pt, outer sep=0pt},
dot/.default = 5pt
}

\def\ch{\operatorname{ch}}
\newcommand{\vtx}[3][]{\node[dot, label={#1}] (#2) at (#3) {};}

\newcommand{\dedge}[3][]{\draw[midarrow] (#2) to[#1] (#3);}

\newcommand\blackslug{\hbox{\hskip 1pt \vrule width 4pt height 8pt depth 1.5pt
        \hskip 1pt}}
\newcommand\bbox{\hfill \quad \blackslug \bigbreak}

\def\DD{\hbox{-}}
\def\CC{\hbox{-}\cdots\hbox{-}}
\def\LL{,\ldots,}
\def\cupcup{\cup\cdots\cup}

\title{Minors of plane digraphs}
\author{Maria Chudnovsky\thanks{Supported by NSF Grants DMS-2348219 and   CCF-2505100,
AFOSR grant FA9550-25-1-0275, and a Guggenheim Fellowship.}\\
Princeton University, Princeton, NJ 08544
\\
\\
Paul Seymour\thanks{Supported by AFOSR grant
FA9550-22-1-0234, and NSF grant  DMS-2154169.}\\
Princeton University, Princeton, NJ 08544}

\date{November 1, 2025; revised \today}

\newtheorem{thm}{}[section]

\newcommand{\Proof}{\noindent{\bf Proof.}\ \ }

\begin{document}
\maketitle
\begin{abstract}

A digraph $H$ is a ``semi-strong minor'' of another, $G$, 
if a subdivision of $H$ can be obtained from a subdigraph of $G$ by contracting strongly-connected subdigraphs to single vertices. We will define a width 
measure of ``plane'' digraphs (that is, drawn in the plane) based on a kind of branch-composition, and show that for every plane digraph $H$, all plane digraphs not containing $H$
as a semi-strong minor have bounded width, while plane digraphs in general have unbounded width. 
\end{abstract}

\section{Introduction}

\begin{figure}[h!]
\centering

\begin{tikzpicture}[scale=.5,auto=left]

\tikzstyle{every node}=[inner sep=1.5pt, fill=black,circle,draw]
\vtx{v11}{1,1};
\node (v12) at (1,2) {};
\node (v13) at (1,3) {};
\node (v14) at (1,4) {};
\node (v21) at (2,1) {};
\node (v22) at (2,2) {};
\node (v23) at (2,3) {};
\node (v24) at (2,4) {};
\node (v31) at (3,1) {};
\node (v32) at (3,2) {};
\node (v33) at (3,3) {};
\node (v34) at (3,4) {};
\node (v41) at (4,1) {};
\node (v42) at (4,2) {};
\node (v43) at (4,3) {};
\node (v44) at (4,4) {};

\draw (v11) -- (v12) -- (v13)-- (v14) -- (v24)--(v23) -- (v22) -- (v21)--(v31) -- (v32) -- (v33)--(v34) -- (v44) -- (v43) -- (v42)--(v41);
\draw (v11) -- (v21);
\draw (v31) -- (v41);
\draw (v12) -- (v22) -- (v32) -- (v42);
\draw (v13) -- (v23) -- (v33) -- (v43);
\draw (v24) -- (v34);

\end{tikzpicture}
\caption{A $4\times 4$ grid.} \label{fig:4grid}
\end{figure}
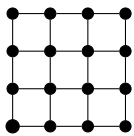

There is a theorem of Robertson and Seymour~\cite{GM5}, that:
\begin{thm}\label{GM5}
For every $k\ge 0$ there exists $\ell\ge 0$ such that every graph with tree-width at least $\ell$ contains a $k\times k$ grid as a
minor. Conversely, for every $\ell\ge 0$ there exists $k\ge 0$ such that every graph that contains a $k\times k$ grid as a
minor has tree-width at least $\ell$.
\end{thm}
(We omit the definitions of ``grid'', ``minor'' and ``tree-width'' since they are well known. See Figure \ref{fig:4grid}.)
See~\cite{diestel} for a simpler proof, and~\cite{chuzhoy} for the best dependence between grid size and tree-width known.

The result \ref{GM5} is equivalent to saying that
for every planar graph $H$, all graphs not containing $H$ as a minor have
bounded tree-width, because every planar graph is a minor of a grid. We are interested in extending this to digraphs.

That raises several issues: what is ``minor containment'' for digraphs? What is ``tree-width'' for digraphs? What is a 
``grid'' for digraphs? But at least we can avoid the first of these issues for the moment.  One can recast the theorem 
of~\cite{GM5} as a statement about the undirected graphs that do not contain a subdivision of a large ``wall'' as a subgraph, and the same idea 
will work for us in digraphs.

\begin{figure}[h!]
\centering

\begin{tikzpicture}[scale=1,auto=left]

\tikzstyle{every node}=[inner sep=1.5pt, fill=black,circle,draw]

\begin{scope}[shift ={(-.2,.2)}]
\node (11) at (1,1) {};
\node (13) at (1,3) {};
\node (15) at (1,5) {};
\node (22) at (2,2) {};
\node (24) at (2,4) {};
\node (26) at (2,6) {};
\node (31) at (3,1) {};
\node (33) at (3,3) {};
\node (35) at (3,5) {};
\node (42) at (4,2) {};
\node (44) at (4,4) {};
\node (46) at (4,6) {};
\node (51) at (5,1) {};
\node (53) at (5,3) {};
\node (55) at (5,5) {};
\node (62) at (6,2) {};
\node (64) at (6,4) {};
\node (66) at (6,6) {};
\end{scope}
\begin{scope}[shift ={(-.2,-.2)}]

\node (12) at (1,2) {};
\node (14) at (1,4) {};
\node (16) at (1,6) {};
\node (21) at (2,1) {};
\node (23) at (2,3) {};
\node (25) at (2,5) {};
\node (32) at (3,2) {};
\node (34) at (3,4) {};
\node (36) at (3,6) {};
\node (41) at (4,1) {};
\node (43) at (4,3) {};
\node (45) at (4,5) {};
\node (52) at (5,2) {};
\node (54) at (5,4) {};
\node (56) at (5,6) {};
\node (61) at (6,1) {};
\node (63) at (6,3) {};
\node (65) at (6,5) {};
\end{scope}

\begin{scope} [shift ={(.2,-.2)}]

\node (119) at (1,1) {};
\node (139) at (1,3) {};
\node (159) at (1,5) {};
\node (229) at (2,2) {};
\node (249) at (2,4) {};
\node (269) at (2,6) {};
\node (319) at (3,1) {};
\node (339) at (3,3) {};
\node (359) at (3,5) {};
\node (429) at (4,2) {};
\node (449) at (4,4) {};
\node (469) at (4,6) {};
\node (519) at (5,1) {};
\node (539) at (5,3) {};
\node (559) at (5,5) {};
\node (629) at (6,2) {};
\node (649) at (6,4) {};
\node (669) at (6,6) {};
\end{scope}

\begin{scope} [shift ={(.2,.2)}]

\node (129) at (1,2) {};
\node (149) at (1,4) {};
\node (169) at (1,6) {};
\node (219) at (2,1) {};
\node (239) at (2,3) {};
\node (259) at (2,5) {};
\node (329) at (3,2) {};
\node (349) at (3,4) {};
\node (369) at (3,6) {};
\node (419) at (4,1) {};
\node (439) at (4,3) {};
\node (459) at (4,5) {};
\node (529) at (5,2) {};
\node (549) at (5,4) {};
\node (569) at (5,6) {};
\node (619) at (6,1) {};
\node (639) at (6,3) {};
\node (659) at (6,5) {};
\end{scope}

\foreach \from/\to in {
119/11,11/12,12/129,129/139,139/13,13/14,14/149, 149/159,159/15,15/16,16/169,
26/269,269/259,259/25,25/24,24/249,249/239,239/23,23/22,22/229,229/219,219/21,
319/31,31/32,32/329,329/339,339/33,33/34,34/349, 349/359,359/35,35/36,36/369,
46/469,469/459,459/45,45/44,44/449,449/439,439/43,43/42,42/429,429/419,419/41,
519/51,51/52,52/529,529/539,539/53,53/54,54/549, 549/559,559/55,55/56,56/569,
66/669,669/659,659/65,65/64,64/649,649/639,639/63,63/62,62/629,629/619,619/61,
61/519,51/419,41/319,31/219,21/119,
129/22,229/32,329/42,429/52,529/62,
63/539,53/439,43/339,33/239,23/139,
149/24,249/34,349/44,449/54,549/64,
65/559,55/459,45/359,35/259,25/159,
169/26,269/36,369/46,469/56,569/66}
\dedge{\from}{\to};

\end{tikzpicture}
\caption{The $6\times 6$ diwall.} \label{fig:diwall}
\end{figure}
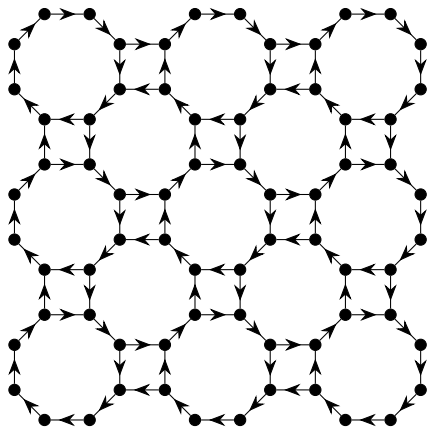

Let us say a digraph $G$ {\em embeds} a digraph $H$ if there is a subdigraph
of $G$ that is isomorphic to a subdivision of $H$.
For $k\ge 2$, even, let us say a {\em $k\times k$ diwall} is a digraph as illustrated in Figure \ref{fig:diwall}. The
subdigraph formed by the horizontal and diagonal edges is the disjoint union of $k$ directed paths, alternating in direction, that we
call the ``horizontal'' paths; and similarly there are
$k$ ``vertical'' paths. Each horizontal
path intersects each vertical path in exactly one edge (the diagonal edges in the figure). If we contract the diagonal edges we obtain a $k\times k$ grid called an {\em alternating grid} that we will discuss later.

Our theorem concerns the structure of digraphs drawn in the 2-sphere that do not embed a large diwall. 

Before we begin with more precise statements and definitions, let us give an idea of what we mean by ``structure'' here.
For undirected graphs, Menger's theorem involves a vertex-cutset of bounded size, but for digraphs, it involves a cutset 
containing a bounded number of vertices and an unbounded number of edges, all crossing the cutset in the wrong direction.
We have something similar. 
For undirected graphs, bounded tree-width means that the graph can be cut up by ``non-crossing'' vertex-cutsets into 
trivial pieces, where the cutsets have bounded size. Our theorem is the same, except that the cutsets now contain both 
a bounded number of vertices and an unbounded number of edges crossing the cutset in prescribed directions. The edges in the cutset
fall into a bounded number of ``intervals'', and all the edges in one 
interval cross the cutset in the same direction. For ``intervals'' to make sense, we assume that 
we are working with cutsets in a digraph drawn in the plane, so 
there is an associated 
circular order on the elements of each cutset, and the intervals are intervals in this order. 
For this reason, our result only gives
information about planar digraphs. We will prove  that for every $k$, every digraph drawn in the plane that does not embed
the 
$k\times k$ diwall can be cut up into trivial pieces by non-crossing cutsets, each consisting of a bounded number of vertices and a bounded number of intervals of edges as discussed. 

The proof is in two stages. First we assume that $G$ is drawn in the plane, and does not embed the 
$k\times k$ diwall, and in addition we have a bound on the ``interleaving'' of $G$, that is, for each vertex $v$, 
the edges incident with $v$ fall into a bounded number of intervals in their cyclic order around $v$, and in each interval all edges enter $v$ or they all leave $v$. In this case, the theorem is simpler, and we do not need to allow the cutsets to contain vertices; we can cut the graph into trivial (that is,
one-vertex) pieces by a non-crossing collection of edge-cutsets, where the edges in each cutset fall into a bounded
number of intervals as before. Then we deduce the stronger theorem, when there is no bound on interleaving.

Later in the paper, we discuss the connection with the well-known digraph grid theorem of Kawarabayashi and 
Kreutzer~\cite{kawa1, kawa2}; and also deduce some consequences for excluding general planar digraphs as minors of large planar digraphs.

\section{Drawings}
We have omitted several definitions, and been sloppy about some others,  in order to state the main result quickly; but now
we need to go back and define things more carefully. So we will basically start again. Graphs and digraphs in this paper
are all finite, and may have parallel edges or loops.
Throughout the paper, $\Sigma$ denotes the 2-sphere, and let us assign an orientation ``clockwise'' to $\Sigma$.
An {\em O-arc} in $\Sigma$ is a subset homeomorphic to
a circle, and a {\em line} in $\Sigma$ is a subset homeomorphic to $[0,1]$. We will not consider drawings in general surfaces; in this paper a {\em drawing} $G$ is 
a graph, such that:
\begin{itemize}
\item  $V(G)\subseteq \Sigma$;
\item  each $e=uv\in E(G)$ is either a line with ends $u\ne v$ in $V(G)$, or an O-arc containing $u=v$;
\item if $v\in V(G)$ and $e\in E(G)$ are not incident in $G$ then $v\notin e$; and
\item if $e,f\in E(G)$ are distinct, every point in $e\cap f$ equals $v$ for some $v\in V(G)$ incident with both $e,f$.
\end{itemize}
We will only consider drawings with no topological pathologies.
A {\em region} is a maximal arc-wise connected subset of $\Sigma$ disjoint from the drawing, and hence
is an open set.
If a drawing $G$ is a digraph we say $G$ is a {\em didrawing}. 
Two drawings or didrawings are {\em homeomorphic} if there is a homeomorphism of $\Sigma$ to itself that takes one to the other.
If $G$ is a digraph, we say a digraph $H$ is a {\em
subdivision} of $G$ if $H$ can be obtained from $G$ by repeatedly replacing an edge $uv$ with a new vertex $w$ and two
edges $uw,wv$
($u,v$ might be equal).  Let us say a digraph $G$ {\em embeds} a digraph $H$ if there is a subdigraph
of $G$ that is isomorphic to a subdivision of $H$.
Similarly,
a didrawing $H$ is a {\em $\Sigma$-subdivision} of a didrawing $G$ if $H$ can be obtained from $G$ by repeatedly replacing an edge $uv$ with a new vertex $w$ and two
edges $uw,wv$, such that $uw\cup wv = uv$.
A didrawing $G$  {\em $\Sigma$-embeds} a didrawing $H$ if there is a subdidrawing of $G$ that is homeomorphic to a 
$\Sigma$-subdivision of $H$.

If $G$ is a digraph, $G^\natural$ denotes the undirected graph underlying $G$.
If $G$ is a drawing, $G^*$ denotes a dual drawing. If $G$ is a didrawing, $G^*$ denotes a didrawing obtained from a dual drawing
${G^\natural}^*$ of $G^\natural$ by directing its edges such that each edge of $G^*$ crosses the corresponding edge of $G$ from left to right (we recall that we fixed some orientation of $\Sigma$, so ``left to right''
makes sense). We call $G^*$ a {\em right-dual didrawing}, and a {\em left-dual didrawing} is defined similarly. (Taking right-dual and then left-dual returns us to the original didrawing.) A {\em dual} drawing is either a right- or left-dual drawing.

\section{The decomposition}

Next, we want to explain the decomposition more precisely. 
The {\em leaves} of a tree are its vertices of degree one, and 
a tree $T$ is {\em ternary} if each of its vertices has degree one or three.

Let $S$ be a finite set with $|S|\ge 2$, let $T$ be a ternary tree,
and let $\tau$ be a bijection from $S$ onto the set of leaves of $T$. We call $(T,\tau)$ a {\em carving} of $S$. 
For each
edge $e$ of $T$ let $T_1(e), T_2(e)$ be the two components of $T\setminus e$. 
We call 
$$(\{v \in  S : \tau(v)\in V(T_1(e))\},\{v \in  S : \tau(v)\in V(T_2(e))\})$$
and its reverse the {\em partitions of $S$ corresponding to $e$ under $(T,\tau)$}.

For instance, if $G$ is a graph, we can ask for a carving $(T,\tau)$ of $E(G)$
such that for each edge $e\in E(T)$, there are at most $k$ vertices incident with edges in both parts
of the partition corresponding to $e$; the smallest $k$ for which such a carving  exists is called the {\em branch-width} of $G$, and was 
introduced in \cite{GM10}. It is within a constant factor of tree-width. In~\cite{ratcatcher}, Seymour and Thomas showed that one can compute the branch-width of a 
planar graph in polynomial time, but it is NP-hard for general graphs. 

Similarly, if $G$ is a graph, we can ask for a carving $(T,\tau)$ of $V(G)$
such that for each edge $e\in E(T)$, there are at most $k$ edges between the two parts
of the partition corresponding to $e$; the smallest $k$ for which such $(T,\tau)$ exists is the {\em carving-width}, 
and again it was shown in~\cite{ratcatcher} that one can compute the carving-width of a 
planar graph in polynomial time, and it is NP-hard for general graphs. The same paper also showed that if $G$ is 2-connected,
there is an ``optimal''  carving $(T,\tau)$ such that for each edge $e\in E(T)$, the set of edges
joining the two parts of the corresponding partition is a {\em bond} of $G$, that is, a minimal edge-cutset of $G$, and so corresponds to a cycle of the dual graph. 

We want something very close to this for didrawings $G$. We will look for carvings $(T,\tau)$ of $V(G)$ where
the corresponding edge-cutsets are bonds of $G$ (that is, of $G^\natural$), but we will measure
the ``size'' of edge-cutsets in a different way.

If $D$ is a cycle of a digraph $G$, some of its vertices might be the head of both incident edges, or the tail of both these edges. We call such a vertex a {\em
change-vertex} of $D$, and the number of change-vertices is the {\em change number} of $D$. Thus, a cycle $D$ has change number zero if and only if it is a
directed cycle; and the change number of every cycle is even. Similarly, if $P$ is a path of a digraph $G$, its {\em change number} is the number of vertices $v$ of $P$ such that $v$ is the head of two edges of $P$ or the tail of two edges of $P$.
If $G$ is a didrawing, the change number of a bond of $G$ is the change number of the corresponding cycle in a dual didrawing.

For $k\ge 1$, we say a digraph $G$ is {\em weakly $k$-connected} or {\em $k$-weak} if $G^\natural$ is $k$-connected.
Let us say a didrawing has {\em  diwidth} at most $k$ if there is a carving $(T,\tau)$ of $V(G)$, such that for each 
edge $e$ of $T$, the set of edges joining the two parts of the 
corresponding partition is a bond of $G$ with change number at most $k$. (Didrawings that are not 2-weak do not admit 
any such ``bond'' carvings, so diwidth is only defined for 2-weak didrawings.)
Then our first theorem says:

\begin{thm}\label{edgecutthm}
For all even $k\ge 2$, and every odd integer $\lambda\ge 1$,  every 2-weak loopless didrawing with interleaving at most $2\lambda$ and with diwidth more than $16\lambda^2k^2$ embeds the $k\times k$ diwall.
\end{thm}

We remark that our interest in diwidth grew from a paper by Berger, Carter and Seymour~\cite{weightable}, where it is shown
that a strongly-connected, 2-weak loopless didrawing is ``weightable'' if and only if it has diwidth at most two. ``Weightable'' means that one can assign
weights to its edges totalling to 1 in every directed cycle. (There are several equivalent conditions, for instance that
a 0/1-weighting exists, or that for every triple of vertices, every directed cycle containing them uses them in the same cyclic order.)

Now let us turn to the decomposition used for our main theorem, when there is no bound on the interleaving. 
Let $G$ be a loopless didrawing. 
A {\em dart} of $G$ is a pair
$(v,e)$ where $v\in V(G)$, $e\in E(G)$ and $e$ is incident with $v$. 
Let $F$ be a simple closed curve in $\Sigma$, such that $F$ has at most one point in common with each
edge, and for each region $r$, either $F\cap r=\emptyset$ or $F\cap r$ is the interior of a line. We call such
a curve $F$ {\em good}. Let $F$ be a good curve, and let $\Delta_1,\Delta_2$ be the two closed discs in $\Sigma$ bounded by $F$.
We say a dart $(v,e)$ {\em belongs to} $\Delta_i$ if either $v$ belongs to the interior of $\Delta_i$, or $v\in F$
and $e\setminus \{v\}$ is a subset of the interior of $\Delta_i$.
For $i = 1,2$, let $D_i$ be the set of darts that belong to $\Delta_i$; we call $(D_1,D_2)$ and $(D_2,D_1)$ the 
{\em dart-partitions from $F$}. A partition $(A,B)$ of the set of all darts is {\em sensible} if $(A,B)$ is a dart-partition from $F$ for some good curve $F$. 

Let $F$ be a good curve, and let $\Delta_1,\Delta_2,D_1,D_2$ be as before. A {\em change-region} for $F$ is a region $r$ incident with two edges $e,f$ that both
cross $F$ (since $F$ is good, there cannot be three edges crossing $F$ and incident with $r$),
such that $\Delta_1$ contains the head of one of $e,f$ and $\Delta_2$ contains the head of the other.
(If $r,e,f$ are as above, it follows that
$F\cap (r\cup e\cup f)$ is a line with ends in the interiors of $e$ and $f$.)
The number of change-regions for $F$ is its {\em change number}.
If $F\cap V(G)=\emptyset$, this coincides with our previous definition of change number for bonds. 
Different good curves $F$ may give rise to the same dart-partition $(D_1,D_2)$, but it follows that all such curves
have the same change number; the number of change-regions for $F$ is determined by
the partition $(D_1,D_2)$. Let us say the {\em cost} of $(D_1,D_2)$ is the sum of the change number of $F$ and 
$|F\cap V(G)|$ (again, it is easy to check that the cost of $F$ is also determined by the dart-partition).

We say a loopless didrawing has {\em dart-width} at most $k$ if there is a carving $(T,\tau)$ of the set of all darts, 
such that for each edge $e$ of $T$, the corresponding dart-partition is sensible and has cost at most $k$.
A digraph is {\em weakly 2-edge-connected} if $G^\natural$ is 2-edge-connected (no didrawing with a cut-edge admits a
sensible dart-partition).
We will prove:
\begin{thm}\label{mainthm}
For all even $k\ge 2$, if $G$ is a loopless, weakly 2-edge-connected didrawing
that does not embed 
the $k\times k$ diwall, then $G$ has dart-width at most $16k^2$.
\end{thm}

If we apply \ref{mainthm} to  
the digraph obtained from an undirected planar graph by replacing each edge by two opposite directed edges, we obtain \ref{GM5} restricted to planar graphs and re-expressed in terms of branch-width.

\section{Patterns of directed cuts}

\begin{figure}[h!]
\centering

\begin{tikzpicture}[scale=.6,auto=left]

\draw [rounded corners](0,-3) rectangle (2,3);
\draw [rounded corners](4,-3) rectangle (6,3);
\draw [rounded corners](8,-3) rectangle (10,3);
\draw [rounded corners](12,-3) rectangle (14,3);
\draw [rounded corners](16,-3) rectangle (18,3);
\draw [rounded corners](20,-3) rectangle (22,3);

\draw[midarrow, line width=1] (1.5,2) to (4.5,1.8);
\draw[midarrow, line width=1] (1.5,.3) to (4.5,-.2);
\draw[midarrow, line width=1] (1.5,-1.8) to (4.5,-2.2);

\draw[midarrow, line width=1] (5.5,1.8) to (8.5,2.3);
\draw[midarrow, line width=1] (5.5,0) to (8.5,0);
\draw[midarrow, line width=1] (5.5,-1.6) to (8.5,-2.3);
\draw[midarrow, line width=1] (5.5,1) to (8.5,.8);

\draw[midarrow, line width=1] (12.5,1.8) to (9.5,2.3);
\draw[midarrow, line width=1] (12.5,0) to (9.5,0);
\draw[midarrow, line width=1] (12.5,-1.6) to (9.5,-1);
\draw[midarrow, line width=1] (12.5,-2.2) to (9.5,-2.5);

\draw[midarrow, line width=1] (13.5,1.8) to (16.5,2.1);
\draw[midarrow, line width=1] (13.5,0) to (16.5,0.2);
\draw[midarrow, line width=1] (13.5,-1.6) to (16.5,-2.3);

\draw[midarrow, line width=1] (20.5,1.8) to (17.5,2.1);
\draw[midarrow, line width=1] (20.5,.30) to (17.5,-.1);
\draw[midarrow, line width=1] (20.5,-1.2) to (17.5,-1.3);
\draw[midarrow, line width=1] (20.5,-2.2) to (17.5,-2.5);

\tikzstyle{every node}=[inner sep=1.5pt, fill=black,circle,draw]
\node (u) at (1,0) {};
\node (v) at (21,0) {};
\tikzstyle{every node}=[]
\draw (u) node [left]           {$u$};
\draw (v) node [right]           {$v$};

\node at (1,.7) {$A_0$};
\node at (5,.7) {$A_1$};
\node at (9,.7) {$A_2$};
\node at (13,.7) {$A_3$};
\node at (17,.7) {$A_4$};
\node at (21,.7) {$A_5$};

\end{tikzpicture}
\caption{A $(u,v)$-multicut with pattern $(1,1,-1,1,-1)$,} \label{fig:multicut}
\end{figure}
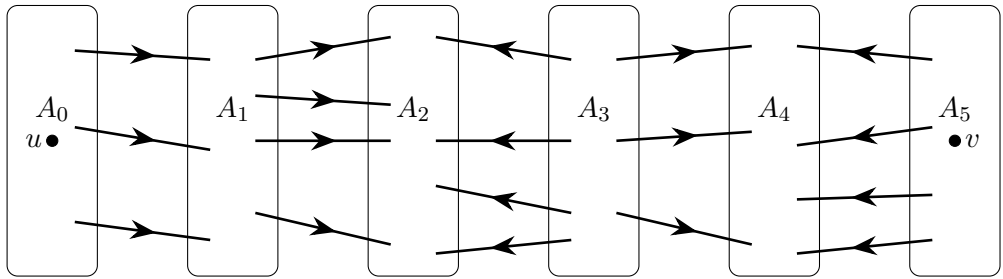

Much of the paper is concerned with proving \ref{edgecutthm}, and we will begin on that now.
A {\em pattern} is a sequence $(\pi_1\LL \pi_{k})$ where $k\ge 1$ and $\pi_i\in \{+1,-1\}$ for $1\le i\le k$.
Let $u,v$ be distinct vertices of a 1-weak digraph $G$. Let us say a {\em $(u,v)$-multicut} is a sequence $(A_0\LL A_k)$ of subsets
of $V(G)$ with $k\ge 1$, such that
\begin{itemize}
\item $A_0\LL A_k$ are
pairwise disjoint and have union $V(G)$;
\item $u\in A_0$ and $v\in A_k$;
\item for $0\le i,j\le k$, if $j\ge i+2$
there are no edges between $A_i$ and $A_j$;
\item for $1\le i\le k$, either there are no edges from $A_{i-1}$ to $A_{i}$ or there are no edges from $A_{i}$ to $A_{i-1}$.
\end{itemize}
(See Figure \ref{fig:multicut}.)

If $(A_0\LL A_k)$ is a $(u,v)$-multicut, let $\pi_1\LL \pi_{k}$ be defined by $\pi_i = 1$ if there are no edges from $A_{i}$ to $A_{i-1}$, and $\pi_i = -1$ if there are no edges from $A_{i-1}$ to $A_{i}$ for $1\le i\le n$. (This is well-defined, since there is some edge between $A_{i-1}, A_{i}$ because $G$ is 1-weak.)
We call the sequence $\pi = (\pi_1\LL \pi_{k})$
the {\em pattern} of the multicut $(A_0\LL A_k)$.

We are interested in when there is a $(u,v)$-multicut with a given pattern.
The main result of this section is:

\begin{thm}\label{getcuts}
Let $u,v$ be distinct vertices of a 1-weak digraph $G$, and let $\pi = (\pi_1\LL \pi_{k})$ be a pattern. Then there is a $(u,v)$-multicut in $G$ with pattern $\pi$ if and only if for every path of $G$ between $u,v$, there is a
$(u,v)$-multicut in $P$ with pattern $\pi$.
\end{thm}
\Proof The ``only if'' part is clear. For the ``if'' part,
we proceed by induction on $k$, and for fixed $k$, by induction on $|E(G)|$. We assume that there is no
$(u,v)$-multicut in $G$ with pattern $\pi$, and we must show that there is a path $P$ between $u,v$ that admits no
$(u,v)$-multicut in $G$ with pattern $\pi$. If $k=1$ the claim is clear, so we assume
that $k\ge 2$. We may assume that there are no loops (by deleting them, since loops make no difference). By reversing all edges of $G$ and negating $\pi$ if necessary, we may assume that $\pi_k=+1$.
\\
\\
(1) {\em If there is an edge of $G$ with
tail $v$, then there is a path $P$ between $u,v$ that admits no
$(u,v)$-multicut in $G$ with pattern $\pi$.}
\\
\\
Suppose that $e=vw \in E(G)$.
If $w=u$, then $e$ makes a path $P$ between $u,v$ that admits no
$(u,v)$-multicut in $G$ with pattern $\pi$ as required. So we assume that $w\ne u,v$. Let $G'$
be obtained from $G$ by contracting $e$, and let $v'$ be the vertex of $G'$ made by identifying $v,w$ under contraction.
Since there is no
$(u,v)$-multicut in $G$ with pattern $\pi$, it follows that there is no
$(u,v')$-multicut in $G'$ with pattern $\pi$. From the inductive hypothesis,
there is a path $P'$ of $G'$ between $u,v'$ that admits no
$(u,v')$-multicut with pattern $\pi$. Hence there is a path $P$ of $G$ between $u,v$, such that either
$E(P) = E(P')$ or $E(P) = E(P')\cup \{e\}$. We claim that in both cases, $P$ satisfies (1).
This is clear if $E(P) = E(P')$, so we assume that  $E(P) = E(P')\cup \{e\}$, and hence
$e$ is the edge of $P$ incident with $v$, and $P\setminus v$ is a path of $G$ between $u,w$. Suppose that $P$ admits a
$(u,v)$-multicut $(A_0\LL A_k)$ with pattern $\pi$. Since $\pi_k=+1$ and $e$ has tail $v$, it follows that 
$w\notin A_{k-1}$; and $w\notin A_i$ for $i\le k-2$ from the definition of a $(u,v)$-multicut. Thus $w\in A_k$, and so contracting
$e$ yields a $(u,v')$-multicut in $P'$ with pattern $\pi$, a contradiction. This proves that $P$ satisfies (1), and so proves (1).

\bigskip

From (1) we may assume that every edge of $G$ incident with $v$ has head $v$. Let $N$ be the set of all vertices adjacent
to $v$, and let $G'$ be obtained from $G$ by contracting all edges between $N,v$. Let $v'$ be the vertex of $G'$ made by
identifying $N\cup \{v\}$ under contraction. Let $\pi'$ be the pattern $(\pi_1\LL \pi_{k-1})$.
\\
\\
(2) {\em There is no
$(u,v')$-multicut in $G'$ with pattern $\pi'$.}
\\
\\
Suppose that
$(A_0'\LL A_{k-1}')$ is a $(u,v')$-multicut in $G'$ with pattern $\pi'$.  Thus $v'\in A_{k-1}'$. Let $A_i = A_i'$ for $0\le i\le k-2$, let $A_{k-1} = (A_{k-1}'\setminus \{v'\})\cup N$,
and let $A_k = \{v\}$. We claim that $(A_0\LL A_{k})$ is a $(u,v)$-multicut in $G$ with pattern $\pi$. To show this, we must check that
\begin{itemize}
\item $A_0\LL A_k$ are pairwise disjoint and have union $V(G)$, and $u\in A_0$ and $v\in A_k$ (and this is all clear);
\item there are no edges between $A_i, A_{k-1}$ with $i\le k-3$ (this is true since there are no edges of $G'$ between 
$A_i', A_{k-1}'$);
\item there are no edges between $A_i, A_{k}$ with $i\le k-2$ (this is true since all edges of $G$ incident with $v$ are between $v,N$);
\item every edge between $A_{k-1}, A_k$ is from $A_{k-1}$ to $A_k$ (this is true since no edge has tail $v$).
\end{itemize}
This shows that $(A_0\LL A_{k})$ is a $(u,v)$-multicut in $G$ with pattern $\pi$, a contradiction, and so proves (2).

\bigskip

From the inductive hypothesis on $k$,
there is a path $P'$ of  $G'$
between $u,v'$ that admits no
$(u,v')$-multicut with pattern $\pi'$. Hence there is a path $Q$ of $G$ between $u, N\cup \{v\}$ with $E(Q) = E(P')$; and
since every vertex in $N$ is an in-neighbour of $v$, there is a path $P$ of $G$ between $u,v$ such that either $E(P) =E(P')$, or $E(P) = E(P')\cup \{e\}$ from some edge $e = qv$, where $q\in N$ is the end of $Q$ in $N$. If $E(P) = E(P')$
then $P$ admits no
$(u,v)$-multicut with pattern $\pi'$, and hence none with pattern $\pi$, as required. Thus, we assume that
$E(P) = E(P')\cup \{e\}$ from some edge $e = qv$. Suppose that $(A_0\LL A_{k})$ is a $(u,v)$-multicut in $P$ with pattern
$\pi$. Hence $(A_0\LL A_{k-2}, A_{k-1}\cup A_k)$ is a $(u,v)$-multicut in $P$ with pattern
$\pi'$. The edge $e$ has an end in $A_k$, and so $q\notin A_0\LL A_{k-2}$. Consequently, both its ends are in $A_{k-1}\cup A_k$;
and therefore contracting $e$ converts  $(A_0\LL A_{k-2}, A_{k-1}\cup A_k)$ to a $(u,v')$-multicut in $P'$
with pattern $\pi'$, a contradiction. This proves that there is no $(u,v)$-multicut in $P$ with pattern
$\pi$, and therefore completes the proof of \ref{getcuts}.~\bbox

Here is a companion lemma:
\begin{thm}\label{makebonds}
Let $u,v$ be distinct vertices of a 1-weak digraph $G$, and suppose there is a $(u,v)$-multicut in $G$ with some pattern $\pi=(\pi_1\LL \pi_k)$.
Then there is a $(u,v)$-multicut $(A_0\LL A_k)$ in $G$ with pattern $\pi$, such that for $1\le i\le k$, $G[A_0\cup A_1\cupcup A_{i-1}]$
and $G[A_i\cupcup A_k]$ are both weakly connected.
\end{thm}
\Proof
Choose a $(u,v)$-multicut $(A_0\LL A_k)$ in $G$ with pattern $\pi$, such that the sum (over $1\le i\le k$) of the number of edges between $A_{i-1}, A_i$
is minimum. Suppose that for some $i\in \{1\LL k\}$, $G[A_0\cup A_1\cupcup A_{i-1}]$ is not connected, and let $X,Y$ be a partition of
$A_0\cup A_1\cupcup A_{i-1}$ into two nonempty subsets such that there are no edges between $X,Y$, and $u\in X$.
Then
$$(A_0\cap X, A_1\cap X\LL A_{i-1}\cap X, A_i\cup Y, A_{i+1}\LL A_k)$$
is also a $(u,v)$-multicut $(A_0\LL A_k)$ in $G$ with pattern $\pi$; so, from the minimality in the choice of $(A_0\LL A_k)$,
there is no edge between $Y, A_i$, contradicting that $G$ is 1-weak. Thus
$G[A_0\cup A_1\cupcup A_{i-1}]$ is connected, and similarly so is  $G[A_i\cupcup A_k]$, for $1\le i\le k$. This proves \ref{makebonds}.\bbox

Thus, \ref{getcuts} reduces the general question to the same question when $G$ is a path.
Let us look more closely at the path case. We are mostly interested in when
the desired sequence $(\pi_1\LL \pi_k)$ is {\em alternating}, that is, $\pi_i, \pi_{i+1}$ have opposite sign for $1\le i\le k-1$.

Let $P$ be a path between $u,v$ in a digraph $G$, and let $u = q_0,q_1\LL q_t = v$ be vertices of $P$ in order (not necessarily all distinct). For $1\le i\le t$ let $P_i$ be the subpath between $q_{i-1}, q_i$, and suppose that each $P_i$ is a directed path.
We say that $(P_1\LL P_t)$ is a {\em concatenation} for $(P, u,v)$.
Now let $\pi = (\pi_1\LL \pi_t)$ be an alternating  pattern; we say that
the concatenation $(P_1\LL P_t)$ is a {\em $\pi$-concatenation} for $(P,u,v)$ if for $1\le i\le t$,
$P_i$ is a directed path from $q_{i-1}$ to $q_i$ if and only if $\pi_i=1$. (A concatenation may be a $\pi$-concatenation for more than one 
pattern $\pi$ since we permit some of the paths $P_i$ to have length zero.) 
If $\pi=(\pi_1\LL \pi_t)$ is a pattern, $-\pi$ denotes the pattern $(-\pi_1\LL -\pi_t)$.

We claim that:
\begin{thm}\label{alternatingcase}
Let $P$ be a path between $u,v$ and let $\pi = (\pi_1\LL \pi_k)$ be an alternating pattern. Then exactly one of the following holds:
\begin{itemize}
\item $P$ admits a $(u,v)$-multicut with pattern $\pi$;
\item there is a $(-\pi)$-concatenation for $(P,u,v)$.
\end{itemize}
\end{thm}
We omit the proof, which is clear.

If $\pi =(\pi_1\LL \pi_k)$ is an alternating pattern of length $k$, let $\pi^{++}$ denote the alternating pattern obtained by adding a new first and last term, that is,
$$(-\pi_1, \pi_1\LL \pi_k,-\pi_k).$$
Note that $-\pi$ is an initial subsequence of $\pi^{++}$.

We also need the following:
\begin{thm}\label{intersect}
Let $G$ be a 1-weak digraph, let $u,v\in V(G)$ be distinct, let $\pi = (\pi_1\LL \pi_k)$ be an alternating pattern, and suppose that
$(A_0\LL A_k)$ is a $\{u,v\}$-multicut
with pattern $\pi$ in $G$. Let $P$ be a path between $u,v$ such that there is a $\pi^{++}$-concatenation for $(P,u,v)$.
Then $P$ contains exactly one edge between $A_{i-1}, A_i$ for $1\le i\le k$.
\end{thm}
\Proof
For $1\le i\le k$, there is an edge of $P$
between $A_{i-1}, A_i$; let $e_i, f_i$ be the first and last such edges in $P$.  An edge $e=ab$ of $P$ is {\em forward} in $(P,u,v)$
if $u,a$ belong to the same component of $P\setminus e$, and {\em backward} otherwise.
\\
\\
(1) {\em for $1\le i\le k$, $e_i, f_i$ are forward edges of $P$ if and only if $\pi_i = 1$.}
\\
\\
The component of $P\setminus e_i$ containing $u$ is a path with no edges between $A_{i-1}, A_i$, and since it has one end $u$, all its vertices
belong to $A_0\cupcup A_{i-1}$. Hence $e_i$ is a forward edge of $P$ if and only if all edges of $G$ between $A_{i-1}, A_i$
are from $A_{i-1}$ to $A_i$, that is, if and only if $\pi_i=1$. This proves the claim for $e_i$, and the claim for $f_i$ is proved similarly, using the component of $P\setminus f_i$ containing $v$. This proves (1),
\\
\\
(2) {\em $e_i = f_i$ for $1\le i\le k$.}
\\
\\
Define $\pi_0 = -\pi_1$, and $\pi_{k+1} = -\pi_k$. Thus $\pi^{++} = (\pi_0\LL \pi_{k+1})$.
Now $e_1\LL e_k$ are all distinct and in order in $P$. Let $(P_0\LL P_{k+1})$ be a $\pi^{++}$-concatenation for $(P,u,v)$.
It follows from (1)
that $e_1\notin E(P_0)$, and therefore $e_1\LL e_k\in E(P_1)\LL E(P_{k+1})$. Similarly $e_2\notin E(P_1)$, so
$e_2\LL e_k\in E(P_2)\LL E(P_{k+1})$; and in general, $e_i\in E(P_{i})\LL E(P_{k+1})$. Similarly, $f_i\in E(P_0)\LL E(P_{i})$;
and since $e_i$ equals or is earlier than $f_i$ in $P$, it follows that $e_i,f_i\in E(P_i)$. Suppose that $e_i\ne f_i$ for some $i$.
If $\pi_i=1$, then there is a directed path of $P_i$ from the head of $e_i$ to the tail of $f_i$, which is impossible since the head of $e_i$ is in $A_i$ and the tail of $f_i$ is in $A_{i-1}$, and all edges of $G$ between $A_{i-1}, A_i$ are from $A_{i-1}$ to $A_i$; and
the other case (when $\pi_i=-1$) is similar and we omit it. This proves (2).

\bigskip

From (2), only one edge of $P$ is between $A_{i-1}, A_i$ for $1\le i\le k$. This proves \ref{intersect}.~\bbox

We deduce:

\begin{thm}\label{minimalcuts}
Let $G$ be a 1-weak digraph, let $u,v\in V(G)$ be distinct, let $\pi = (\pi_1\LL \pi_k)$ be an alternating pattern, and suppose that
$(A_0\LL A_k)$ is a $\{u,v\}$-multicut
with pattern $\pi$ in $G$. Let $f\in E(G)$, and suppose that there is no $\{u,v\}$-multicut
with pattern $\pi$ in $G/f$ ($G/f$ denotes the digraph obtained by contracting $f$). Then for some $i$, $f$ is between $A_{i-1}$ and $A_i$.
Moreover, there is a path $P$ between $u,v$ and a $\pi^{++}$-concatenation for $(P,u,v)$, such that $f$ is the only edge of $P$
between $A_{i-1}, A_i$.
\end{thm}
\Proof
By \ref{getcuts}, there is a path $P'$ in $G/f$ between $u,v$ such that $(P',u,v)$ admits a 
$(-\pi)$-concatenation $(P_1\LL P_k)$. This no longer exists in $G$,
so there is a path $P$ of $G$ between $u,v$ with $f\in E(P)$, such that $P'$ is obtained from $P$ by contracting $f$. Hence
$(P,u,v)$ admits a $\pi^{++}$-concatenation; and hence the claim follows from \ref{intersect}. This proves \ref{minimalcuts}.~\bbox

\section{Box systems}

For $k\ge 2$, even, let us say a {\em $k \times k$ diwall layout}
in a digraph $G$ consists of two sequences $P_1\LL P_k$ and $Q_1\LL Q_k$ of paths of $G$, with the following properties:
\begin{itemize}
\item $P_1\LL P_k$ are vertex-disjoint directed paths of $G$, and so are $Q_1\LL Q_k$;
\item for all $i,j$ with $1\le i,j\le k$, $P_i\cap Q_j$ is a path $R_{i,j}$ say;
\item if $i$ is odd, the paths $R_{i,1}, R_{i,2}\LL R_{i,k}$ are in order in $P_i$, and if $i$ is even they are in reverse order;
\item if $j$ is even, the paths $R_{1,j}, R_{2,j}\LL R_{k,j}$ are in order in $Q_j$, and if $j$ is odd they are in reverse order.
\end{itemize}
Thus,  if $P_1\LL P_k$ and $Q_1\LL Q_k$ form a $k\times k$ diwall layout then their union embeds a $k\times k$ diwall; and
conversely, if $G$ embeds a $k\times k$ diwall then it contains a $k\times k$ diwall layout.

\begin{figure}[h!]
\centering

\begin{tikzpicture}[scale=.8,auto=left]

\tikzstyle{every node}=[inner sep=1.5pt, fill=black,circle,draw]
\node (11) at (1,1) {};
\node (12) at (1,2) {};
\node (13) at (1,3) {};
\node (14) at (1,4) {};
\node (21) at (2,1) {};
\node (22) at (2,2) {};
\node (23) at (2,3) {};
\node (24) at (2,4) {};
\node (31) at (3,1) {};
\node (32) at (3,2) {};
\node (33) at (3,3) {};
\node (34) at (3,4) {};
\node (41) at (4,1) {};
\node (42) at (4,2) {};
\node (43) at (4,3) {};
\node (44) at (4,4) {};
\node (15) at (1,5) {};
\node (16) at (1,6) {};
\node (25) at (2,5) {};
\node (26) at (2,6) {};
\node (35) at (3,5) {};
\node (36) at (3,6) {};
\node (45) at (4,5) {};
\node (46) at (4,6) {};
\node (51) at (5,1) {};
\node (11) at (1,1) {};
\node (52) at (5,2) {};
\node (53) at (5,3) {};
\node (54) at (5,4) {};
\node (55) at (5,5) {};
\node (56) at (5,6) {};
\node (61) at (6,1) {};
\node (62) at (6,2) {};
\node (63) at (6,3) {};
\node (64) at (6,4) {};
\node (65) at (6,5) {};
\node (66) at (6,6) {};

\foreach \from/\to in {
11/12,12/13,13/14,14/15,15/16,
21/22,22/23,23/24,24/25,25/26,
31/32,32/33,33/34,34/35,35/36,
46/45,45/44,44/43,43/42,42/41,
56/55,55/54,54/53,53/52,52/51,
66/65,65/64,64/63,63/62,62/61,
61/51,51/41,41/31,31/21,21/11,
62/52,52/42,42/32,32/22,22/12,
63/53,53/43,43/33,33/23,23/13,
14/24,24/34,34/44,44/54,54/64,
15/25,25/35,35/45,45/55,55/65,
16/26,26/36,36/46,46/56,56/66}
\dedge{\from}{\to};

\begin{scope}[shift ={(9,0)}]
\node (11) at (1,1) {};
\node (12) at (1,2) {};
\node (13) at (1,3) {};
\node (14) at (1,4) {};
\node (21) at (2,1) {};
\node (22) at (2,2) {};
\node (23) at (2,3) {};
\node (24) at (2,4) {};
\node (31) at (3,1) {};
\node (32) at (3,2) {};
\node (33) at (3,3) {};
\node (34) at (3,4) {};
\node (41) at (4,1) {};
\node (42) at (4,2) {};
\node (43) at (4,3) {};
\node (44) at (4,4) {};
\node (15) at (1,5) {};
\node (16) at (1,6) {};
\node (25) at (2,5) {};
\node (26) at (2,6) {};
\node (35) at (3,5) {};
\node (36) at (3,6) {};
\node (45) at (4,5) {};
\node (46) at (4,6) {};
\node (51) at (5,1) {};
\node (11) at (1,1) {};
\node (52) at (5,2) {};
\node (53) at (5,3) {};
\node (54) at (5,4) {};
\node (55) at (5,5) {};
\node (56) at (5,6) {};
\node (61) at (6,1) {};
\node (62) at (6,2) {};
\node (63) at (6,3) {};
\node (64) at (6,4) {};
\node (65) at (6,5) {};
\node (66) at (6,6) {};

\foreach \from/\to in {
11/12,12/13,13/14,14/15,15/16,
26/25,25/24,24/23,23/22,22/21,
31/32,32/33,33/34,34/35,35/36,
46/45,45/44,44/43,43/42,42/41,
51/52,52/53,53/54,54/55,55/56,
66/65,65/64,64/63,63/62,62/61,
61/51,51/41,41/31,31/21,21/11,
12/22,22/32,32/42,42/52,52/62,
63/53,53/43,43/33,33/23,23/13,
14/24,24/34,34/44,44/54,54/64,
65/55,55/45,45/35,35/25,25/15,
16/26,26/36,36/46,46/56,56/66}
\dedge{\from}{\to};
\end{scope}

\end{tikzpicture}
\caption{The $6\times 6$ semi-grid and alternating grid.} \label{fig:semigrid}
\end{figure}
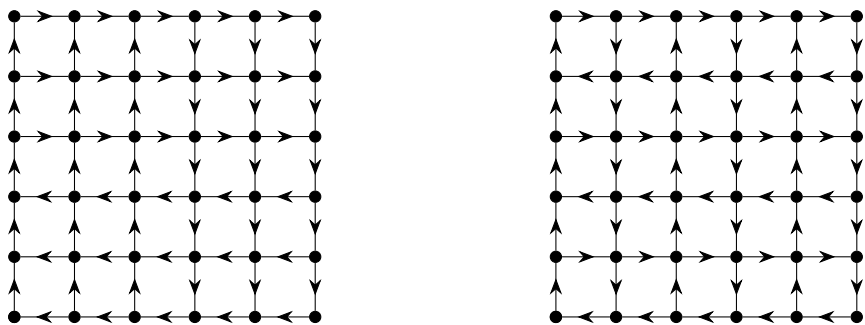

We need the {\em $k_1\times k_2$ alternating grid}, defined for $k_1,k_2\ge 2$, even, and illustrated in Figure \ref{fig:semigrid} (which also shows the ``semi-grid'', that we do not need yet). It has $k_1$ horizontal and $k_2$ vertical paths; 
they 
alternate in direction, and the first vertex of the top horizontal
path is the last vertex of a vertical path (this last condition is just to make it unique). We denotes the
$k_1\times k_2$ alternating grid by $\Gamma_{k_1,k_2}$.

Third, for $k_1,k_2\ge 2$, even, let us say a {\em $k_1\times k_2$ box system} in a digraph $G$ is a map $\eta$ from the 
set of vertices and edges of $H=\Gamma_{k_1,k_2}$ with the following properties:
\begin{itemize}
\item for each $v\in V(H)$, $\eta(v)$ is a subset of $V(G)$, and if $u,v\in V(H)$ are distinct then $\eta(u)\cap \eta(v) = \emptyset$;
\item for each $e=uv\in E(H)$, $\eta(e)\in E(G)$, and the tail and head of $\eta(e)$ are in $\eta(u), \eta(v)$ respectively;
\item for each $v\in V(H)$, if $e,f\in E(H)$ have head $v$ and tail $v$ respectively, then there is a directed path in $G[\eta(v)]$ from the head of $\eta(e)$ to the tail of $\eta(f)$.
\end{itemize}
Evidently a $k \times k$ box system is some kind of grid-like object in $G$. One might hope to extract a $k \times k$ diwall layout
from it, by choosing for each horizontal path $P_i$ of $\Gamma_{k,k}$ a directed path $P_i'$ in $G$ containing all the edges 
$\eta(e) \:(e\in E(P_i))$,
joining these edges via appropriate paths in $G[\eta(v)]$ for $v\in V(P_i)$, and defining $Q_j'$
for each vertical path $Q_j$ of $\Gamma_{k,k}$ similarly. But this does not work, because we cannot arrange that $P_i\cap Q_j$ is a path for each $i,j$.
Nevertheless, we can extract $k \times k$ diwall layout from a $k\times 3k$ box system, as we show next.

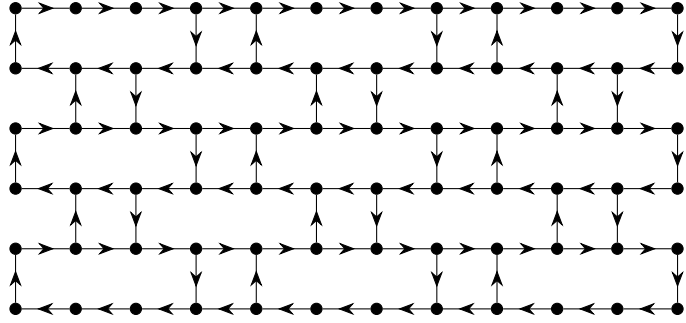
\begin{figure}[h!]
\centering

\begin{tikzpicture}[scale=.8,auto=left,rotate=270]

\tikzstyle{every node}=[inner sep=1.5pt, fill=black,circle,draw]
\node (11) at (1,1) {};
\node (12) at (1,2) {};
\node (13) at (1,3) {};
\node (14) at (1,4) {};
\node (15) at (1,5) {};
\node (16) at (1,6) {};
\node (17) at (1,7) {};
\node (18) at (1,8) {};
\node (19) at (1,9) {};
\node (110) at (1,10) {};
\node (111) at (1,11) {};
\node (112) at (1,12) {};

\node (21) at (2,1) {};
\node (22) at (2,2) {};
\node (23) at (2,3) {};
\node (24) at (2,4) {};
\node (25) at (2,5) {};
\node (26) at (2,6) {};
\node (27) at (2,7) {};
\node (28) at (2,8) {};
\node (29) at (2,9) {};
\node (210) at (2,10){};
\node (211) at (2,11) {};
\node (212) at (2,12) {};

\node (31) at (3,1) {};
\node (32) at (3,2) {};
\node (33) at (3,3) {};
\node (34) at (3,4) {};
\node (35) at (3,5) {};
\node (36) at (3,6) {};
\node (37) at (3,7) {};
\node (38) at (3,8) {};
\node (39) at (3,9) {};
\node (310) at (3,10){};
\node (311) at (3,11) {};
\node (312) at (3,12) {};

\node (41) at  (4,1) {};
\node (42) at  (4,2) {};
\node (43) at  (4,3) {};
\node (44) at  (4,4) {};
\node (45) at  (4,5) {};
\node (46) at  (4,6) {};
\node (47) at  (4,7) {};
\node (48) at  (4,8) {};
\node (49) at  (4,9) {};
\node (410) at (4,10){};
\node (411) at (4,11) {};
\node (412) at (4,12) {};

\node (51) at  (5,1) {};
\node (52) at  (5,2) {};
\node (53) at  (5,3) {};
\node (54) at  (5,4) {};
\node (55) at  (5,5) {};
\node (56) at  (5,6) {};
\node (57) at  (5,7) {};
\node (58) at  (5,8) {};
\node (59) at  (5,9) {};
\node (510) at (5,10){};
\node (511) at (5,11) {};
\node (512) at (5,12) {};

\node (61) at  (6,1) {};
\node (62) at  (6,2) {};
\node (63) at  (6,3) {};
\node (64) at  (6,4) {};
\node (65) at  (6,5) {};
\node (66) at  (6,6) {};
\node (67) at  (6,7) {};
\node (68) at  (6,8) {};
\node (69) at  (6,9) {};
\node (610) at (6,10){};
\node (611) at (6,11) {};
\node (612) at (6,12) {};

\foreach \from/\to in {
11/12,12/13,13/14,14/15,15/16,16/17,17/18,18/19,19/110, 110/111,111/112,
212/211,211/210,210/29,29/28,28/27,27/26,26/25,25/24,24/23,23/22,22/21,
31/32,32/33,33/34,34/35,35/36,36/37,37/38,38/39,39/310,310/311,311/312,
412/411,411/410,410/49,49/48,48/47,47/46,46/45,45/44,44/43,43/42,42/41,
51/52,52/53,53/54,54/55,55/56,56/57,57/58,58/59,59/510,510/511,511/512,
612/611,611/610,610/69,69/68,68/67,67/66,66/65,65/64,64/63,63/62,62/61,
61/51,52/42,41/31,32/22,21/11,
65/55,56/46,45/35,36/26,25/15,69/59,510/410,
49/39,310/210,29/19,
14/24,23/33,34/44,43/53,54/64,
18/28,27/37,38/48,47/57,58/68,
112/212,211/311,312/412,411/511,512/612}
\dedge{\from}{\to};

\end{tikzpicture}
\caption{The $6\times 6$ diwall, redrawn. If we subdivide once the odd edges of the horizontal paths we obtain a subdigraph of the $k\times 3k$ alternating grid.} \label{fig:6diwallagain}
\end{figure}

Figure \ref{fig:6diwallagain} gives another picture of the $k\times k$ diwall, redrawn without diagonal edges. Each of
the ``vertical'' paths now consists of the edges in two consecutive columns of the picture, together with edges of the horizontal paths
joining them.

\begin{thm}\label{boxestodiwall}
Let $k\ge 2$, even, and suppose that there is a $k\times 3k$ box system $\eta$ in $G$; then there is a $k\times k$ diwall 
layout in $G$, and hence a $k\times k$ diwall in $G$.
\end{thm}
\Proof
Let $P_1\LL P_k$ and $Q_1\LL Q_{3k}$ be the horizontal and vertical paths of $\Gamma_{k,3k}$, numbered in order, and
for $1\le i\le k$ and $1\le j\le 3k$ let $v_{ij}$ be the vertex in both $P_i, Q_j$.
For every pair $(e,f)$ where $e,f\in E(\Gamma_{k,3k})$ and the head of $e$ equals the tail of $f$, say $v$, let
$L(e,f)$ be a directed path of $G$ with first edge $\eta(e)$, last edge $\eta(f)$, and all internal vertices in $\eta(v)$,
chosen subject to the following rule:
\begin{itemize}
\item if $v$ is incident with two edges in a horizontal path of $\Gamma_{k,3k}$, say $x,y$, and $|\{e,f\}\cap \{x,y\}| = 1$, then $L(e,f)$ is chosen such that $L(e,f)\cap L(x,y)$ is a path.
\end{itemize}
(In other words, $L(e,f)$ consists of a subpath of $L(x,y)$ together with a minimal path between $L(x,y)$ and whichever of $e,f$ is not in $\{x,y\}$.)

For $1\le i\le k$, let the edges of $P_i$ be $e_1\LL e_{3k-1}$ in order; then 
$$L(e_1,e_2)\cup L(e_2,e_3)\cupcup L(e_{3k-1},e_{3k})$$
is a directed path of $G$ that we call $P_i'$, for $1\le i\le k$. 

For $i=0$ modulo 6, there is a directed path ($R_i$, say) of $\Gamma_{k,3k}$ with vertices 
$$v_{1,i}\DD v_{2,i}\DD  v_{2,i-1}\DD  v_{2,i-2}\DD v_{3,i-2}\DD  v_{3,i-1}\DD  v_{3,i}\DD v_{4,i}\DD v_{4,i-1}\DD v_{4, i-2}\DD v_{5,i-2}\CC v_{k,i}$$
in order, 
and for $i = 1$ mod 6, there is a directed path $R_i$ with vertices
$$v_{k,i}\DD v_{k-1,i}\DD v_{k-1,i+1}\DD v_{k-1,i+2}\DD v_{k-2,i+2}\DD v_{k-2,i+1}\DD v_{k-2,i}\DD v_{k-3,i}\CC v_{1,i}.$$
For $1\le j\le 3k$ with $j\in \{0,1\}$ modulo 6, let the edges of $R_j'$ be $e_1\LL e_{3k-5}$ in order, and let 
$R_j'$ be the path
$$L(e_1,e_2)\cup L(e_2,e_3)\cupcup L(e_{3k-6},e_{3k-5}).$$
Then $P_1'\LL P_k'$ and $R_4', R_5', R_{10}', R_{11}'\LL R_{3k-1}'$ form a $k\times k$ diwall layout. (The drawing in
Figure \ref{fig:6diwallagain} might make this easier to see.) This proves \ref{boxestodiwall}.~\bbox


\section{Paths across a disc}

We need several results about paths across  a didrawing in a disc, so let us first set up some notation.
If $G$ is a didrawing, two cycles $C_1,C_2$ of $G$ are {\em non-crossing} if $C_i$ bounds two closed discs 
$\Delta_i,\Delta_i'$ in $\Sigma$ for $i = 1,2$, and one of $\Delta_1,\Delta_1'$ includes one of $\Delta_2,\Delta_2'$.
Let $G$ be a didrawing, and let $v\in V(G)$, such that no loop is incident with $v$.  
By an {\em alternating ring} through $v$, we mean a sequence $(C_1\LL C_k)$ of directed cycles 
of $G$, each containing $v$, pairwise edge-disjoint and pairwise non-crossing, such that, if $e_i, f_i$ are the edges of $C_i$
with tail $v$ and head $v$ respectively, then 
$$f_1,e_2,f_3,e_4\LL e_k, f_k, e_{k-1}\LL f_2,e_1$$ 
(if $k$ is even) or
$$f_1,e_2,f_3,e_4\LL f_k, e_k, f_{k-1}\LL f_2,e_1$$ 
(if $k$ is odd) are in anticlockwise cyclic order around $v$. We call $k$ its {\em size}. 
\begin{figure}[H]
\centering

\begin{tikzpicture}[scale=1.7,auto=left]
\path[use as bounding box] (-1.5,-1) rectangle (1.5,2);

\tikzstyle{every node}=[inner sep=1.5pt, fill=black,circle,draw]
\node (v) at (0,0) {};
\node (f1) at  ({cos(-70)}, {sin(-70)}) {};
\node (e2) at  ({cos(-30)}, {sin(-30)}) {};
\node (f3) at  ({cos(30)}, {sin(30)}) {};
\node (e4) at  ({cos(70)}, {sin(70)}) {};
\node (f4) at  ({cos(110)}, {sin(110)}) {};
\node (e3) at  ({cos(150)}, {sin(150)}) {};
\node (f2) at  ({cos(210)}, {sin(210)}) {};
\node (e1) at  ({cos(250)}, {sin(250)}) {};

\draw[midarrow, line width=1] (v) to (e1);
\draw[midarrow, line width=1] (v) to (e2);
\draw[midarrow, line width=1] (v) to (e3);
\draw[midarrow, line width=1] (v) to (e4);
\draw[midarrow, line width=1] (f1) to (v);
\draw[midarrow, line width=1] (f2) to (v);
\draw[midarrow, line width=1] (f3) to (v);
\draw[midarrow, line width=1] (f4) to (v);

\draw[midarrow, dashed, out=150, in = 30, looseness = 25] (e1) to (f1);
\draw[midarrow, dashed,  out=35, in=145, looseness = 7] (e2) to (f2);
\draw[midarrow, dashed,  out=110, in=70, looseness = 1.3] (e3) to (f3);
\draw[midarrow, dashed, out=160, in = 20, looseness = .9] (e4) to (f4);
PP
\draw[thick, dotted] ({.5*cos(-27)}, {.5*sin(-27)}) arc [start angle=-27, end angle=27, radius=.5];
\draw[thick, dotted] ({.5*cos(153)}, {.5*sin(153)}) arc [start angle=153, end angle=207, radius=.5];

\tikzstyle{every node}=[];
\node at (-.35,-.6) {$e_1$};
\node at (.35,-.6) {$f_1$};
\node at (-.7,-.23) {$f_2$};
\node at (.7,-.26) {$e_2$};
\node at (-.7,.15) {$e_{k-1}$};
\node at (.75,.13) {$f_{k-1}$};
\node at (-.4,.7) {$f_k$};
\node at (.42,.7) {$e_k$};

\draw[dotted,thick] (0,1.25) to (0,1.6);

\end{tikzpicture}
\caption{An alternating ring. There may be more edges incident with $v$ that are not part of the ring.} \label{fig:altring}
\end{figure}
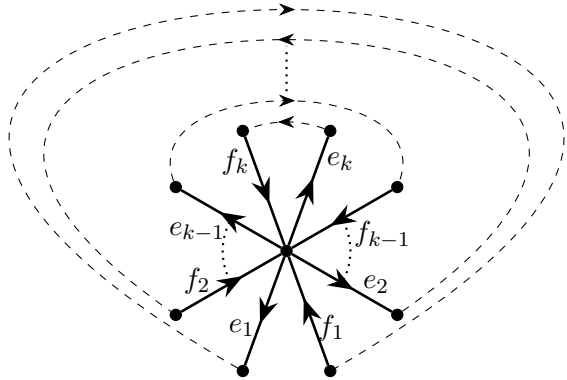

Let $G$ be a 1-weak didrawing, and let $v\in V(G)$, such that no loop is incident with $v$, and such that $G\setminus v$
is 1-weak. Let $G^*$ be a dual didrawing, and
for each vertex, edge or region $x$ of $G$, let $\phi(x)$ be the corresponding region, edge or vertex of $G^*$.
Consequently, the edges of $G$ incident with $v$ are mapped by $\phi$ to the edges of a cycle $C$ of $G^*$ that bounds the region $\phi(v)$ of $G^*$. 
Let $R_1,R_2,R_3,R_4$ be four subpaths of $C$, pairwise edge-disjoint and with $E(R_1\cup R_2\cup R_3\cup R_4) = E(C)$, numbered in 
anticlockwise circular order in $C$ around $\phi(v)$.  An alternating ring $(C_1\LL C_k)$ in $G$ through $v$ is {\em from $R_2$ to $R_4$} if, with 
$e_i, f_i\; (1\le i\le k)$ defined as before, 
\begin{itemize}
\item if $i$ is odd then $\phi(e_i)\in E(R_4)$ and $\phi(f_i)\in E(R_2)$;
\item if $i$ is even then $\phi(e_i)\in E(R_2)$ and $\phi(f_i)\in E(R_4)$.
\end{itemize}
\begin{figure}[h!]
\centering

\begin{tikzpicture}[scale=1,auto=left]
\draw [rounded corners](-2,-1.9) rectangle (2,1.9);

\draw[midarrow, line width=1] (-1.6,1.5) to (-1.6, 2.5);
\draw[midarrow, line width=1] (-.7, 2.5) to (-.7,1.5);
\draw[midarrow, line width=1] (.7,1.5) to (.7, 2.5);
\draw[midarrow, line width=1] (1.6, 2.5) to (1.6,1.5);

\draw[midarrow, line width=1]  (-1.6, -2.5) to (-1.6,-1.5);
\draw[midarrow, line width=1] (-.7,-1.5) to (-.7, -2.5);
\draw[midarrow, line width=1] (.7, -2.5) to (.7,-1.5);
\draw[midarrow, line width=1] (1.6,-1.5) to (1.6, -2.5);

\node at (-1.9,2.2) {$f_1$};
\node at (-1,2.2) {$e_2$};
\node at (1.1,2.2) {$f_{k-1}$};
\node at (1.9,2.2) {$e_k$};
\node at (1.8,-2.3) {$f_k$};
\node at (1.1,-2.3) {$e_{k-1}$};
\node at (-1,-2.3) {$f_2$};
\node at (-1.9,-2.3) {$e_1$};

\node at (0, 1.6) {$R_2$};
\node at (-2.3, 0) {$R_1$};
\node at (0, -1.6) {$R_4$};
\node at (2.3, 0) {$R_3$};

\draw[dotted,thick] (-.5,2.2) to (.5,2.2);
\draw[dotted,thick] (-.5,-2.2) to (.5,-2.2);

\draw[midarrow, dashed] (-1.6, -1.5) to (-1.6, 1.5);
\draw[midarrow, dashed] (-.7, 1.5) to (-.7, -1.5);
\draw[midarrow, dashed] (.7, -1.5) to (.7, 1.5);
\draw[midarrow, dashed] (1.6, 1.5) to (1.6, -1.5);

\end{tikzpicture}
\caption{The end-edges of the members of an alternating ring from $R_2$ to $R_4$. 
The labelled edges are now in clockwise order, since they are all incident with $v$, which is outside the cycle $C$ and not drawn. The dashed paths are edge-disjoint and non-crossing.} \label{fig:altringagain}
\end{figure}
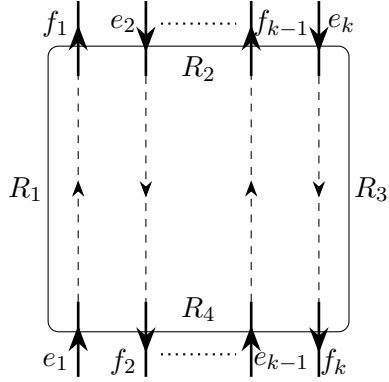

\begin{thm}\label{getalternatingring}
With notation as above, let $k\ge 1$ be an integer. Then either:
\begin{itemize}
\item there is a path $P$ of $G^*$ between $V(R_1)$ and $V(R_3)$
with change number less than $k$; or
\item there is an alternating ring $(C_1\LL C_k)$ in $G$ through $v$ from $R_2$ to $R_4$.
\end{itemize}
\end{thm}
\Proof
Let $\pi = (\pi_1\LL \pi_k)$ be the alternating pattern with $k$ terms and with $\pi_1=1$. 
We assume the statement of the first bullet is false. Consequently, in the digraph $H$ obtained from $G^*$ by contracting
the edges of $R_1$ (making a vertex $u$) and contracting the edges of $R_3$ (making a vertex $w$), there is no path $P$ between
$u,w$ such that $(P,u,w)$ admits a $\pi$-concatenation. By \ref{getcuts} and \ref{makebonds},
there is a $(u,w)$-multicut $(A_0\LL A_k)$ in $H$ with pattern $\pi$, such that for $1\le i\le k$, $H[A_0\cup A_1\cupcup A_{i-1}]$
and $H[A_i\cupcup A_k]$ are both connected. Hence, the set of edges of $H$ between $A_{i-1}$ and $A_i$ is a bond 
of $H$, and so its preimage under $\phi$ is the set of edges of a directed cycle $C_i$ of $H^*$ (and hence of $G$) 
containing $v$.
These cycles are pairwise edge-disjoint, and pairwise non-crossing (since each is the boundary of a set of regions of $H$, 
and these sets are nested by inclusion)
and so $(C_1\LL C_k)$ is an alternating ring in $G$. This proves \ref{getalternatingring}.~\bbox

An alternating ring $(C_1\LL C_k)$ through $v$ is {\em disjointed} if $C_1\LL C_{k_1}$ are pairwise vertex-disjoint except for $v$.
We deduce:

\begin{thm}\label{disjointring}
With notation as before, let $k\ge 1$ be an integer, and suppose that $G$ has interleaving at most $2\lambda$, where $\lambda$ is an odd integer. Then either:
\begin{itemize}
\item there is a path $P$ of $G^*$ between $V(R_1)$ and $V(R_3)$
with change number less than $\lambda k$; or
\item there is a disjointed alternating ring $(C_1\LL C_k)$ in $G$ through $v$ from $R_2$ to $R_4$.
\end{itemize}
\end{thm}
\Proof
We assume the first bullet is false. 
From \ref{getalternatingring}, there is an alternating ring $(B_1\LL B_{k\lambda})$ through $v$
and from $R_2$ to $R_4$. We claim that for $1\le i<j\le \lambda k$, if $j-i\ge \lambda$ then $V(B_i\cap B_j) = \{v\}$;
because suppose that some vertex $u\ne v$ belongs to $V(B_i\cap B_j)$. Then $u$ also belongs to $B_{i'}$ for all $i'$ with $i\le i'\le j$, and so $G$
has interleaving at least $2(j-i+1)$, a contradiction. Let $C_i = B_{(i-1)\lambda+1}$ for $1\le i\le k$; then
$C_1\LL C_{k}$ are pairwise vertex-disjoint except for $v$, and since $\lambda$ is odd, it follows that
$(C_1\LL C_{k})$ is an alternating ring through $v$ from $R_2$ to $R_4$. This proves \ref{disjointring}.~\bbox

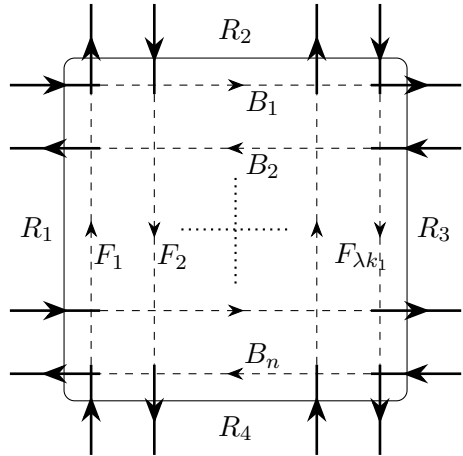
\begin{figure}[h!]
\centering

\begin{tikzpicture}[scale=1.2,auto=left]
\draw [rounded corners](-1.9,-1.9) rectangle (1.9,1.9);

\draw[midarrow, line width=1] (-1.6,1.5) to (-1.6, 2.5);
\draw[midarrow, line width=1] (-.9, 2.5) to (-.9,1.5);
\draw[midarrow, line width=1] (.9,1.5) to (.9, 2.5);
\draw[midarrow, line width=1] (1.6, 2.5) to (1.6,1.5);

\draw[midarrow, line width=1]  (-1.6, -2.5) to (-1.6,-1.5);
\draw[midarrow, line width=1] (-.9,-1.5) to (-.9, -2.5);
\draw[midarrow, line width=1] (.9, -2.5) to (.9,-1.5);
\draw[midarrow, line width=1] (1.6,-1.5) to (1.6, -2.5);

\node at (0, 2.2) {$R_2$};
\node at (-2.2, 0) {$R_1$};
\node at (0, -2.2) {$R_4$};
\node at (2.2, 0) {$R_3$};

\draw[midarrow, dashed] (-1.6, -1.5) to (-1.6, 1.5);
\draw[midarrow, dashed] (-.9, 1.5) to (-.9, -1.5);
\draw[midarrow, dashed] (.9, -1.5) to (.9, 1.5);
\draw[midarrow, dashed] (1.6, 1.5) to (1.6, -1.5);

\draw[midarrow, line width=1] (-1.5,-1.6) to (-2.5,-1.6);
\draw[midarrow, line width=1] (-2.5,-.9) to (-1.5,-.9);
\draw[midarrow, line width=1] (-1.5, .9) to (-2.5,.9);
\draw[midarrow, line width=1] (-2.5, 1.6) to (-1.5,1.6);

\draw[midarrow, line width=1]  (2.5, -1.6) to (1.5, -1.6);
\draw[midarrow, line width=1] (1.5,-.9) to (2.5,-.9);
\draw[midarrow, line width=1] (2.5,.9) to (1.5,.9);
\draw[midarrow, line width=1] (1.5, 1.6) to (2.5,1.6);

\draw[dotted,thick] (0,-.6) to (0,.6);
\draw[dotted,thick] (-.6,0) to (.6,0);

\draw[midarrow, dashed] (1.5,-1.6) to (-1.5,-1.6);
\draw[midarrow, dashed] (-1.5, -.9) to (1.5,-.9);
\draw[midarrow, dashed] (1.5,.9) to (-1.5,.9);
\draw[midarrow, dashed] (-1.5, 1.6) to (1.5, 1.6);

\node at (.3, 1.4) {$B_1$};
\node at (.3, .7) {$B_2$};
\node at (.3, -1.4) {$B_n$};
\node at (-1.4, -.3) {$F_1$};
\node at (-.7, -.3) {$F_2$};
\node at (1.4, -.3) {$F_{\lambda k_1}$};

\end{tikzpicture}
\caption{The $B_i$'s are vertex-disjoint, and the $F_j$'s are edge-disjoint and non-crossing. The paths are in the orders
shown, but the intersections might be much more complicated than indicated. }
\label{fig:disjointpaths}
\end{figure}

We deduce:
\begin{thm}\label{getdisjointpaths}
With notation as before, suppose also that the didrawing $G$ has interleaving at most $2\lambda$ where $\lambda$ is odd.
Let $k_1,k_2\ge 1$ be integers, and let $n=2\lambda k_1k_2+3k_2-2\lambda k_1-2$.
Suppose that every path of $G^*$ between 
$V(R_1)$ and $V(R_3)$ has change number at least $\lambda k_1$, and there is a disjointed alternating ring $(B_1\LL B_n)$
through $v$ 
from $R_3$ to $R_1$.
Then there is a disjointed alternating ring $(C_1\LL C_{k_1})$ 
through $v$ from $R_2$ to $R_4$, and a disjointed alternating ring $(D_1\LL D_{k_2})$ through $v$ from $R_3$ to $R_1$,
such that 
for $1\le h\le k_1$ and $1\le i<j\le k_2$, there is a subpath of $C_i$ with one end $v$ that contains all vertices of 
$C_h\cap D_i$ and no vertices of $C_h\cap D_j$ except $v$. 
\end{thm}
\Proof We proceed by induction on $|E(G)|$.
By \ref{getalternatingring} there is an alternating ring  of size $\lambda k_1$ through $v$ and from $R_2$ to $R_4$, say $(F_1\LL F_{\lambda k_1})$. (See Figure \ref{fig:disjointpaths}.) 
We may assume:
\\
\\
(1) {\em Let $1\le h\le \lambda k_1$, and let $1\le j,j'\le n$ with $j'-j> 2(\lambda k_1+1)$.
Then there is a subpath of $F_h$ with one end $v$ that contains all vertices of
$F_h\cap B_j$ and no vertices of $F_h\cap B_{j'}$ except $v$.}
\\
\\
Since $(B_1\LL B_{n})$ is a disjointed alternating ring, and the path $F_h\setminus v$ contains a vertex of 
$B_{j+\lambda k_1+1}$ and a vertex of  $B_{j+\lambda k_1+2}$, there is a subpath of $F_h\setminus v$ between 
$V(B_{j+\lambda k_1+1})$ and $V(B_{j+\lambda k_1+2})$ with no internal vertex in either of these sets. Let $e$ be an edge of this subpath. Thus $e$ belongs to none of $B_1\LL B_{\lambda k_1k_2}$. 
From the inductive hypothesis
on $|E(G)|$, we may assume that there is no alternating ring in $G\setminus e$ of size $\lambda k_1$ through $v$ and from $R_2$ to $R_4$. Hence,
from \ref{minimalcuts}, there is a path $P$ of $G^*$ between $V(R_1),V(R_3)$ with change number at most $\lambda k_1$   
such that $e$ is the only edge in $E(F_h)$ mapped into $E(P)$ by $\phi$. Consequently, there is a partition $(A,B)$ of 
$V(G)\setminus \{v\}$, such that the edges of $G\setminus \{v\}$ between $A,B$ are precisely the edges of $G\setminus \{v\}$ 
mapped by $\phi$ into $E(P)$.
In particular, $e$ is the only edge of $F_h$ between $A,B$. Hence $F_h[A\cup \{v\}$ is a path $Q_1$ say, and 
$F_h[B\cup \{v\}]$ is a path $Q_2$. From the choice of $e$, each of $B_1\LL B_{j+\lambda k_1+1}$ has a vertex in 
$V(Q_1)\setminus {v}$, and each of $B_{j+\lambda k_1+2}\LL B_{\lambda k_1k_2})$ has a vertex in $V(Q_2)\setminus \{v\}$. 

Let $I$ be the set of all $i\in \{1\LL \lambda k_1k_2\}$ such that some edge of $B_i$ is between $A,B$.
It follows 
from planarity that
$I$ is an interval (that is, if $i_1\le i_2\le i_3$ and $i_1,i_3\in I$ then $i_2\in I$). Since $P$ has change number at most
$\lambda k_1$, it follows that $|I|\le \lambda k_1+1$. If $j\in I$ then, since the path $P$ intersects each of 
$D_j\LL D_{j+\lambda k_1+1}$, it follows that $j,j+1\LL j+\lambda k_1+1\in I$, contradicting that $|I|\le \lambda k_1+1$. Thus
$j\notin I$. Similarly, if $j'\in I$, then $j',j'-1\LL j+\lambda k_1+2\in I$, contradicting that $|I|\le \lambda k_1+1$ 
since $j'-(j+\lambda k_1+2)+1> \lambda k_1+1$; and hence $j'\notin I$. Since no edge of $B_j$ is between $A,B$, and
$B_j$ has a vertex in $V(Q_1)\setminus {v}\subseteq A$, it follows that $V(B_j)\subseteq A\cup \{v\}$, and in particular
all vertices of $F_h\cap B_j$ belong to $Q_1$. Similarly, all vertices of $F_h\cap B_{j'}$ belong to $Q_2$ and therefore
do not belong to $Q_1$. This proves (1).

\bigskip

Let $\mu = 2(\lambda k_1+1)+1$, and for $1\le i\le k_2$ let $D_i = B_{(i-1)\mu+1}$ (this is well-defined, from the choice of $n$.) 
It follows that for $1\le h\le \lambda k_1$ and $1\le i<j\le k_2$, there is a subpath $P$ of $F_h$ with one end $v$ that contains all vertices of
$F_h\cap D_i$ and no vertices of $F_h\cap D_j$ except $v$. For $1\le i\le k_1$ let $C_i = F_{(i-1)\lambda +1}$;
then $(C_1\LL C_{k_1})$ and $(D_1\LL D_{k_2})$ are disjointed alternating rings. This proves \ref{getdisjointpaths}.~\bbox

\begin{thm}\label{ringtowall}
With notation as before,
let $k\ge 2$ be even,  and 
suppose that there is a disjointed alternating ring $(C_1\LL C_{k})$
through $v$ from $R_2$ to $R_4$, and a disjointed alternating ring $(D_1\LL D_{3k})$ through $v$ from $R_3$ to $R_1$,
such that
for $1\le h\le k$ and $1\le i<j\le 3k$, there is a subpath $P$ of $C_i$ with one end $v$ that contains all vertices of
$C_h\cap D_i$ and no vertices of $C_h\cap D_j$ except $v$. Then $G$ embeds the
$k\times k$ diwall.
\end{thm}
\Proof This is immediate from \ref{boxestodiwall}. Alternatively, 
for all $i\in \{1\LL k-1\}$ and all $j\in \{1\LL 3k\}$, let $L_{i,j}$ be a minimal subpath of $D_j\setminus \{v\}$ 
between $V(C_i), V(C_{i+1})$.  Then the paths $C_1\setminus \{v\}\LL C_k\setminus \{v\}$ and the paths
$L_{i,j}$ where $1\le i\le k$ and $1\le j\le 3k$ and either
\begin{itemize}
\item $i$ is odd and $j$ is congruent to 0 or 1 modulo 6, or
\item $i$ is even and $j$ is congruent to 3 or 4 modulo 6
\end{itemize}
form a $k\times k$ diwall layout.
This proves \ref{ringtowall}.~\bbox

Putting these pieces together, we deduce:

\begin{thm}\label{getwall}
Let $k\ge 2$ be even. 
Let $G$ be a 1-weak didrawing, and let $v\in V(G)$, such that no loop is incident with $v$, and such that $G\setminus v$
is 1-weak.  Suppose that $G$ has interleaving at most $2\lambda$, where $\lambda$ is an odd integer. Let $G^*$ be a dual drawing, and
for each vertex, edge or region $x$ of $G$, let $\phi(x)$ be the corresponding region, edge or vertex of $G^*$.
Let $C$ be the cycle of $G^*$ with edge set all edges $\phi(e)$ where $e\in E(G)$ is incident with $v$ in $G$. Let
$n= 4\lambda k^2-6\lambda k +6k-6$,
and suppose that:
\begin{itemize}
\item $R_1,R_2,R_3,R_4$ are edge-disjoint subpaths of $C$, numbered in cyclic order in $C$;
\item every path of $G^*$ between $V(R_1)$ and $V(R_3)$ has change number at least $\lambda n$,
\item every path of $G^*$ between $V(R_2)$ and $V(R_4)$ has change number at least $\lambda k$,
\end{itemize}
Then $G$ embeds the $k\times k$ diwall.
\end{thm}
\Proof By \ref{disjointring}, there is a disjointed alternating ring $(B_1\LL B_n)$
through $v$
from $R_3$ to $R_1$. By \ref{getdisjointpaths}, 
there is a disjointed alternating ring $(C_1\LL C_{k})$
through $v$ from $R_2$ to $R_4$, and a disjointed alternating ring $(D_1\LL D_{3k})$ through $v$ from $R_3$ to $R_1$,
such that
for $1\le h\le k$ and $1\le i<j\le 3k$, there is a subpath of $C_i$ with one end $v$ that contains all vertices of
$C_h\cap D_i$ and no vertices of $C_h\cap D_j$ except $v$.
By \ref{ringtowall}, $G$ embeds the 
$k\times k$ diwall. This proves \ref{getwall}.~\bbox

\section{Shortcuts across a cycle}

\begin{thm}\label{pushcut}
Let $G$ be a 2-weak loopless didrawing, and let $C$ be a cycle of $G$ with change number at most $n$ where $n$ is even. 
Suppose that $G^*$ has interleaving at most $2\lambda$, where $n\ge 4k+2\lambda+10$, and let $\Delta$ be a disc bounded by $C$. 
Suppose that for every path $P$ of $G[\Delta]$ with distinct ends in $V(C)$ and with no internal vertex or edge in $C$, 
there is a subpath $Q$ of $C$ between the ends of $P$ such that $P\cup Q$ has change number more than $n$.
Then either:
\begin{itemize}
\item the interior of $\Delta$ is a region of $G$; or
\item there are four subpaths $R_1,R_2,R_3,R_4$ of $G$, pairwise edge-disjoint, with union $C$, and numbered in order in $C$, such that there is no path in $G$ between $V(R_1), V(R_3)$ with change number at most $k$, and there 
is no path in $G$ between $V(R_2), V(R_4)$ with change number at most $k$.
\end{itemize}
\end{thm}
\Proof We assume the first bullet is false. For a path or cycle $P$ we define $\ch(P)$ to be the change number of $P$.
\\
\\
(1) {\em $\ch(C)\ge n-2\lambda$.}
\\
\\
Each edge of $C$ is incident with a region of $G$ included in $\Delta$. If all these regions are equal, then 
the interior of $\Delta$ is a region of $G$ (because $G$ is 2-weak and loopless) and the first bullet holds, a 
contradiction. Thus  
they are not all equal; and hence there is a subpath $P$ of $C$ and a region $r\subseteq \Delta$ of $G$, such that 
some edge of $C$ is incident with $r$, and every vertex or edge of $C$ incident with $r$ belongs to $V(P)$.
Choose such $P,r$
with $P$ minimal. It follows that every edge of $P$ is incident with $r$ (because if some edge of $P$ is incident with a 
different region $r'\subseteq \Delta$, then all vertices and edges of $C$ incident with $r'$ belong to $P$, from planarity, and for the same reason, not both ends of $P$ are incident with both $r,r'$; this contradicts
the minimality of $P$). The boundary of $r$ is a cycle, since $G$ is 2-weak and loopless, and hence consists of the union of 
$P$ and a path $Q$ of $G[\Delta]$ joining the ends of $P$ with no internal vertex or edge in $C$. 

Let $P'\ne P$ be the second path of $C$ between the ends of $P$. 
Now 
$\ch(P\cup Q)$ is at most $2\lambda\le n$, from the hypothesis, and so $\ch(P'\cup Q)>n$
(and hence at least $n+2$). But $\ch(P'\cup Q)\le \ch(P')+\ch(Q)+2$, and so $\ch(C)\ge \ch(P')\ge n-2\lambda$.
This proves (1). 

\bigskip

From (1), there are four subpaths $R_1\LL R_4$ of $C$, pairwise edge-disjoint and with union $C$, numbered in cyclic order in $C$, and each with change number at least $(n-2\lambda)/4 -2>k$. 
We may assume (for a contradiction) that $R_1\LL R_4$ do not satisfy the theorem;
let $P$ be a path in $G[\Delta]$, either between $V(R_2), V(R_4)$ or between $V(R_1), V(R_3)$,
with 
\begin{itemize}
\item change number as small as possible (and hence at most $k$);  
\item subject to that, with $C\cup P$ minimal; and
\item subject to that, with length as small as possible. 
\end{itemize}
We may assume that $P$ is between $V(R_1), V(R_3)$ from 
the symmetry. From the minimality of the length of $P$, the only vertex of $P$ in $V(R_1)$ is one of its ends, say $r_1$, and its only vertex in $V(R_3)$ is its other end 
$r_3$. Moreover, not both $V(R_2), V(R_4)$ contain internal vertices of $P$, again from the minimality of the length of $P$, so we
may assume that every internal vertex of $P$ in $V(C)$ belongs to $R_2$. Let the components of $P\cap C$ be 
$Q_0\LL Q_m$, numbered in order in $P$, where $r_1\in V(Q_0)$ and $r_3\in V(Q_m)$. Thus, $m\ge 1$, and either $Q_0$ has only one
vertex, or $r_1$ is the common end of $R_1,R_2$, and the same for $Q_m$. For $1\le i\le m$, let $P_i$ be the 
the minimal subpath of $P$ between $V(Q_{i-1}), V(Q_i)$. Thus, $Q_0, P_1, Q_1, P_2\LL P_m, Q_m$ are edge-disjoint subpaths
of $P$ with union $P$, numbered in order in $P$, and $P_1\LL P_m$ each have at least one edge. For $1\le i\le m$, let
$S_i$ be the subpath of the path $R_1\cup R_2\cup R_3$ between the ends of $P_i$, and let $S_i'$ be the other subpath of $C$ between these two vertices. Thus, 
$Q_0, S_1,Q_1\LL S_m, Q_m$ are edge-disjoint subpaths
of $R_1\cup R_2\cup R_3$ with union a path between $R_1,R_3$, numbered in order in this path, and $S_1\LL S_m$ each have at least one edge. 

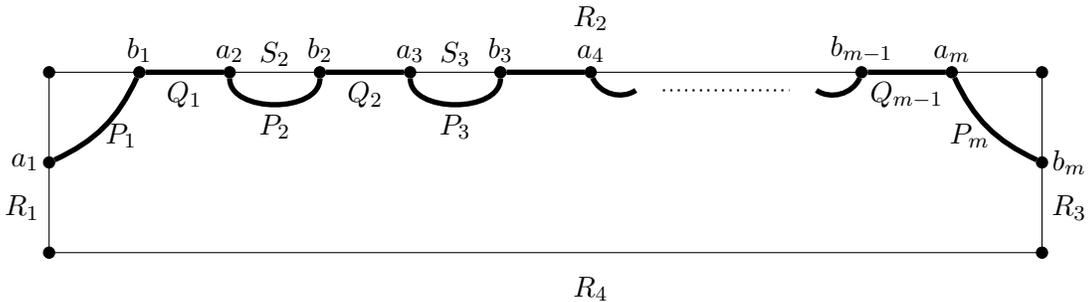
\begin{figure}[h!]
\centering

\begin{tikzpicture}[scale=1.2,auto=left]
\draw (0,0) rectangle (11,2);
\node at (6, 2.6) {$R_2$};
\node at (-.3, .5) {$R_1$};
\node at (6, -.4) {$R_4$};
\node at (11.3, .5) {$R_3$};

\tikzstyle{every node}=[inner sep=1.5pt, fill=black,circle,draw]

\node (z1) at (0,2) {};
\node (z2) at (11,2) {};
\node (z3) at (11,0) {};
\node (z4) at (0,0) {};

\node (a1) at (0,1) {};
\node (b1) at (1,2) {};
\node (a2) at (2,2) {};
\node (b2) at (3,2) {};
\node (a3) at (4,2) {};
\node (b3) at (5,2) {};
\node (a4) at (6,2) {};
\node (b4) at (9,2) {};
\node (a5) at (10,2) {};
\node (b5) at (11,1) {};

\draw[line width=2] (a1) to [bend right = 20] (b1);
\draw[line width=2] (a2) to [bend right = 90] (b2);
\draw[line width=2] (a3) to [bend right = 90] (b3);
\draw[line width=2] (a4) to [bend right = 50] (6.5,1.8);

\draw[line width=2] (8.5, 1.8) to [bend right = 50] (b4);
\draw[line width=2] (a5) to [bend right = 20] (b5);

\draw[line width=2] (b1) to (a2);
\draw[line width=2] (b2) to (a3);
\draw[line width=2] (b3) to (a4);
\draw[line width=2] (b4) to (a5);

\draw[dotted, thick] (6.8,1.8) to (8.2, 1.8);

\tikzstyle{every node}=[]
\draw (a1) node [left] {$a_1$};
\draw (b1) node[above] {$b_1$};
\draw (a2) node[above] {$a_2$};
\draw (b2) node[above] {$b_2$};
\draw (a3) node[above] {$a_3$};
\draw (b3) node[above] {$b_3$};
\draw (a4) node[above] {$a_4$};
\draw (b4) node[above] {$b_{m-1}$};
\draw (a5) node[above] {$a_m$};
\draw (b5) node[right] {$b_m$};

\node at (.8,1.3) {$P_1$};
\node at (2.5, 1.4) {$P_2$};
\node at (4.5, 1.4) {$P_3$};
\node at (10.2, 1.3) {$P_m$};
\node at (1.5, 1.75) {$Q_1$};
\node at (3.5, 1.75) {$Q_2$};
\node at (9.5, 1.75) {$Q_{m-1}$};

\node at (2.5, 2.2) {$S_2$};
\node at (4.5, 2.2) {$S_3$};

\end{tikzpicture}
\caption{The path $P$ and various derived paths. The picture might be misleading, because $P_1$ might have both ends in $R_2$, and then $Q_0$ might have edges; and the same for $P_m, Q_m$. }
\label{fig:jumps}
\end{figure}

For $1\le i\le m$, no internal vertex or edge of $P_i$ belongs to $C$, and so from the hypothesis, one of 
$\ch(P_i\cup S_i), \ch(P_i\cup S_i')>n$.
Since $S_i\cup S_i'=C$ has change number at most $n$, 
and the
change numbers of cycles are even, the change number
of $P_i$ is at least that of one of $S_i, S_i'$. But $\ch(S_i')\ge (n-2\lambda)/4 -2$
since $S_i'$ includes $R_4$, and since $\ch(P_i)\le k< n-2\lambda)/4 -2$, it follows that
$\ch(P_i)\ge \ch(S_i)$. (Consequently the latter is at most $k$, and so $R_2\not\subseteq S_i$, and therefore
$m\ge 2$.) 

If $v$ is a vertex of a path or cycle $T$, we say $\delta(v,T) = 1$ if $v$ is a change-vertex of $T$, and $\delta(v,T) = 0$
otherwise. 
For $1\le i\le n$, let the ends of $P_i$ be $a_i, b_i$, where $a_i$ is between $r_1, b_i$ in $P$.
(Thus $a_i\ne b_i$, but possibly $b_i = a_{i+1}$ if $i<m$.)
\\
\\
(2) {\em For $1\le i\le m$, 
$$\ch(P_i)+\delta(a_i, P_i\cup S_i') + \delta(b_i, P_i\cup S_i')\ge \ch(S_i) + \delta(a_i, C) + \delta(b_i, C)+2.$$}
\noindent
As we have seen, $\ch(P_i\cup S_i')\ge \ch(C)+2$; but 
$$\ch(P_i\cup S_i) = \ch(P_i)+\delta(a_i, P_i\cup S_i') + \delta(b_i, P_i\cup S_i')+\ch(S_i)$$
and 
$$\ch(C)=\ch(S_i) + \delta(a_i, C) + \delta(b_i, C)+\ch(S_i'),$$
and the claim follows. This proves (2).
\\
\\
(3) {\em $Q_0$ has no edges, and $Q_m$ has no edges.}
\\
\\
Suppose that $Q_0$ has an edge. Then $r_1$ is the common end of $R_1,R_2$, and $Q_0$ is a subpath of $R_2$ between $r_1$ 
and $a_1$. Since $m\ge 2$,
$a_1$ is an internal vertex of $R_2$.  From the minimality of $C\cup P$, the union of 
$Q_0\cup S_1$ and the subpath of $P$ between $b_1$ and $r_3$ has change number more (and hence at least two more, 
since $E(Q_0)\ne \emptyset$) 
than that of $P$; and since 
the union of $Q_0\cup P_1$ and  the subpath of $P$ between $b_1$ and $r_3$ equals $P$, we deduce that
$$\delta(a_1,C) + \ch(S_1) + \delta(b_1,C)\ge \delta(a_1, P) + \ch(P_1) + \delta(b_1,P_1\cup S_2)+2.$$
But from (1), 
$$\ch(P_1)+\delta(a_1, P_1\cup S_1') + \delta(b_1, P_1\cup S_1')\ge \ch(S_1) + \delta(a_1, C) + \delta(b_1, C)+2.$$
Adding, we deduce:
\begin{align*}
&\delta(a_1,C) + \ch(S_1) + \delta(b_1,C) + \ch(P_1)+\delta(a_1, P_1\cup S_1') + \delta(b_1, P_1\cup S_1')\\
\ge~&\delta(a_1, P) + \ch(P_1) + \delta(b_1,P_1\cup S_2)+2 + \ch(S_1) + \delta(a_1, C) + \delta(b_1, C)+2.
\end{align*}
which simplifies to 
$$ \delta(a_1, P_1\cup S_1') + \delta(b_1, P_1\cup S_1')
\ge \delta(a_1, P) + \delta(b_1,P_1\cup S_2)+4;$$
and this is impossible since $\delta(a_1, P_1\cup S_1'), \delta(b_1, P_1\cup S_1')\le 1$. Thus $Q_0$ has no edges and similarly
$Q_m$ has no edges. This proves (3).

\bigskip

Let $I=\{i\in \{1\LL m-1\}: b_i = a_{i+1}\}$ and $J= \{1\LL m-1\}\setminus I$. 
It follows that 
$$k\ge \ch(P) \ge  \sum_{1\le i\le m}\ch(P_i) + \sum_{i\in J} \ch(Q_i).$$
Moreover, the change number of $S_1\cupcup S_m$ (which includes $R_2$) equals
$$\sum_{1\le i\le m}\ch(S_i) + \sum_{i\in I}\delta(a_i,C) + \sum_{i\in J}(\delta(a_i,C) + \ch(Q_i) + \delta(b_i,C))\ge \ch(R_1)\ge (n-2\lambda)/4 -2>k,$$
and by summing these two inequalities we obtain 
$$\sum_{1\le i\le m}\ch(S_i) + \sum_{i\in I}\delta(a_i,C) + \sum_{i\in J}(\delta(a_i,C) + \delta(b_i,C)) > 
\sum_{1\le i\le m}\ch(P_i).$$
On the other hand, (2) implies that
$$\ch(P_i)\ge \ch(S_i) + \delta(a_i, C) + \delta(b_i, C)$$
for $1\le i\le m$, and 
by summing for $1\le i\le m$ we deduce:
$$\sum_{1\le \le m} \ch(P_i) \ge \sum_{1\le i\le m} \ch(S_i) + 
\sum_{1\le i\le m}(\delta(a_i, C) + \delta(b_i, C)).$$
Combining with the previous inequality, we deduce:
\begin{align*}
&\sum_{1\le i\le m}\ch(S_i) + \sum_{i\in I}\delta(a_i,C) + \sum_{i\in J}(\delta(a_i,C) + \delta(b_i,C))
+\sum_{1\le \le m} \ch(P_i) \\
>&\sum_{1\le i\le m}\ch(P_i) + \sum_{1\le i\le m} \ch(S_i) + 
\sum_{1\le i\le m}(\delta(a_i, C) + \delta(b_i, C)),
\end{align*}
which simplifies to a contradiction.

This proves that $R_1\LL R_4$ satisfy the theorem, and so proves \ref{pushcut}.~\bbox

By combining \ref{pushcut} and \ref{getwall}, we deduce:

\begin{thm}\label{finalgetwall}
Let $k\ge 2$ be even. 
Let $G$ be a 2-weak loopless didrawing with interleaving at most $2\lambda$, where $\lambda$ is an odd integer. 
Let $\ell=16\lambda^2k^2$.
Let $G^*$ be a dual drawing, let $C$ be a cycle of $G^*$ with change number at most $\ell$, and let $\Delta$ be a disc bounded by $C$. 
Suppose that, for every path $P$ of $G^*[\Delta]$ with distinct ends in $V(C)$ and with no internal vertex or edge in $C$,
there is a subpath $Q$ of $C$ between the ends of $P$ such that $P\cup Q$ has change number more than $\ell$.
Then either:
\begin{itemize}
\item the interior of $\Delta$ is a region of $G$; or
\item $G$ has a subdigraph that is a directed subdivision of the $k\times k$ diwall.
\end{itemize}
\end{thm}
\Proof Let $n= 4\lambda k^2-6\lambda k +6k-6$, and $K=\lambda n$. Thus, $\ell \ge  4K+2\lambda+10$. We assume the first bullet is false, and so by \ref{pushcut} applied to $G^*$, 
there are four subpaths $R_1,R_2,R_3,R_4$ of $G$, pairwise edge-disjoint, with union $C$, and numbered in order in $C$, 
such that there is no path in $G^*$ between $V(R_1), V(R_3)$ with change number at most $K$, and there
is no path in $G^*$ between $V(R_2), V(R_4)$ with change number at most $K$. Since $K\ge \lambda \max(n,k)$, 
\ref{getwall} implies that 
$G$ has a subdigraph that is a directed subdivision of the $k\times k$ diwall. This proves \ref{finalgetwall}~\bbox

Now we can deduce \ref{edgecutthm}, which we restate:
\begin{thm}\label{edgecutthm2}
For all even $k\ge 2$, and every odd integer $\lambda\ge 1$,  every loopless 2-weak didrawing with interleaving at 
most $2\lambda$ and with diwidth more than $16\lambda^2k^2$ embeds the $k\times k$ diwall.
\end{thm}
\Proof
Since $G$ is 2-weak, it has a vertex $v$, and the set of edges incident with $v$ is a bond of $G^\natural$. Hence there 
is a ternary tree $T$ and a surjective map $\tau$ from $V(G)$ to the set of leaves of $T$, such that for each edge $e$ of $T$, 
if $T_1,T_2$ are the components of $T\setminus e$, then the set of edges of $G$ between 
$$(\{v \in  V(G) : \tau(v)\in V(T_1)\},\{v \in  V(G) : \tau(v)\in V(T_2)\})$$
is a bond of $G$ with change number at most $16\lambda^2k^2$ (because we can take $T$ to be a two-vertex tree, and map $v$ to one vertex of $T$ and all others to the other vertex of $T$, and $2\lambda\le 16\lambda^2k^2$). Choose $T$ maximal with this property. We claim that $\tau$ is a bijection. Suppose not; then for some leaf $t$ of $T$, at least two vertices of $G$
are mapped to $t$ by $\tau$. Let $e$ be the edge of $T$ incident with $t$, and let $A=\{v \in  V(G) : \tau(v)=t\}$
and $B=V(G)\setminus A$. Thus, $|A|\ge 2$, and the set of edges of $G$ between $A,B$ is a bond of $G$ with change number at most $16\lambda^2k^2$. Let 
$G^*$ be a dual didrawing, and let $C$ be the cycle of $G^*$ corresponding to this bond. Thus $C$ has change number
at most $16\lambda^2k^2$. Let $\Delta$ be the closed disc in $\Sigma$ bounded by $C$ that includes $A$.
Since $|A|\ge 2$, $\Delta$ includes more than one region of $G^*$, and so by 
by \ref{finalgetwall}, we may assume that
there is a path $P$ of $G^*[\Delta]$ with distinct ends in $V(C)$ and with no internal vertex or edge in $C$,
such that $P\cup Q$ has change number at most $\ell$ for each of the subpaths $Q$ of $C$ between the ends of $P$.
Let $Q_1,Q_2$ be these two subpaths, and for $i = 1,2$, let $\Delta_i$ be the closed disc in $\Delta$ bounded by 
$P\cup Q_i$, and let $A_i$ be the set of vertices of $G$ in $\Delta_i$ (that is, that correspond under the duality to a 
region of $G^*$ included in $\Delta_i$). Add two new vertices $t_1,t_2$ to $T$, both adjacent to $t$, forming a ternary tree $T'$, and define $\tau'$ by
\begin{align*}
\tau'(v) &= \tau(v) \text{ if }v\in B\\
\tau'(v) &= t_1 \text{ if } v\in A_1\\
\tau'(v) &= t_2 \text{ if } v\in A_2.
\end{align*}
Since $P\cup Q_i$ is a cycle of $G^*$ with change number at most $16\lambda^2k^2$, it follows that the set of edges between
$A_i$ and $V(G)\setminus A_i$ is a bond of $G$ with change number at most $16\lambda^2k^2$, contrary to the maximality of 
$T$. This proves \ref{edgecutthm2}.~\bbox

\section{The main result}

Now we use \ref{edgecutthm2} to deduce our main result, \ref{mainthm}. We need the following construction.
Let $G$ be a 2-edge-connected loopless 
drawing, and let $J$ be obtained from $G$ as follows. Subdivide each edge $uv$ twice, so it becomes a directed path 
of length three. This makes a 2-edge-connected didrawing $G_1$ say.
Next, for each vertex $v$ of $G$, let $a_1\LL a_n$ be its neighbours in $G_1$ in clockwise order, and add edges to $G_1$
from $a_i$ to $a_{i+1}$ for $1\le i\le n$ (where $a_{n+1} = a_1$). This makes $G_2$. Now delete all the vertices 
of $G$ from $G_2$, forming $J$. Thus, $J$ is a 2-weak didrawing, and each vertex $v$ of $G$ corresponds to a directed 
cycle $C_v$
of $J$ that bounds a region of $J$ in a clockwise direction; and each edge of $G$ corresponds to an edge of $J$ between the corresponding two cycles. Each vertex of $J$ has indegree two and outdegree one, or vice versa. 
We can recover $G$ from $J$ by contracting the edges of the cycles $C_v\;(v\in V(G))$. 

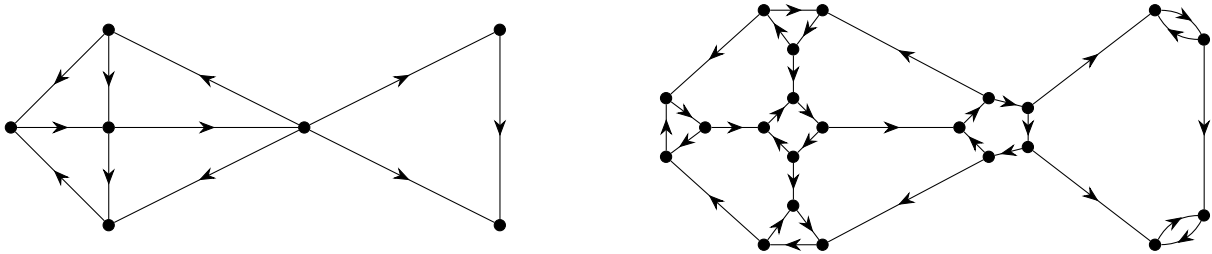
\begin{figure}[h!]
\centering

\begin{tikzpicture}[scale=1.3,auto=left]

\tikzstyle{every node}=[inner sep=1.5pt, fill=black,circle,draw]
\node (3) at (3,0) {};
\node (2) at (2,0) {};
\node (1) at (3,1) {};
\node (4) at (3,-1) {};
\node (5) at (5,0) {};
\node (6) at (7,1) {};
\node (7) at (7,-1) {};

\node (12) at (9.7,1.2) {};
\node (15) at (10.3,1.2) {};
\node (13) at (10,.8) {};
\node (42) at (9.7,-1.2) {};
\node (45) at (10.3,-1.2) {};
\node (43) at (10,-.8) {};
\node (21) at (8.7,0.3) {};
\node (23) at (9.1,0) {};
\node (24) at (8.7, -0.3) {};
\node (31) at (10,.3) {};
\node (32) at (9.7,0) {};
\node (34) at (10,-.3) {};
\node (35) at (10.3,0) {};
\node (51) at (12, .3) {};
\node (53) at (11.7,0) {};
\node (54) at (12,-.3) {};
\node (56) at (12.4,.2) {};
\node (57) at (12.4,-.2) {};
\node (65) at (13.7,1.2) {};
\node (67) at (14.2,.9) {};
\node (75) at (13.7,-1.2) {};
\node (76) at (14.2,-.9) {};
\foreach \from/\to in {
1/2,1/3,5/1,2/3,4/2, 3/4, 3/5,5/4, 5/6,5/7,6/7}
\dedge{\from}{\to};

\foreach \from/\to in {
12/21,13/31,51/15,23/32,42/24, 34/43, 35/53,54/45, 56/65,57/75,67/76}
\dedge{\from}{\to};

\foreach \from/\to in {
12/15,15/13,13/12, 21/23,23/24,24/21, 42/43,43/45,45/42, 32/31,31/35,35/34,34/32, 51/56,56/57,57/54,54/53,53/51}
\dedge{\from}{\to};

\foreach \from/\to in {65/67,67/65,76/75,75/76} \dedge[bend left  = 30]{\from}{\to};

\end{tikzpicture}
\caption{Constructing $J$ from $G$.} \label{fig:blowup}
\end{figure}

The didrawing $J$ is 2-weak and loopless, and has interleaving at most two; so by \ref{edgecutthm2} applied to it, we deduce that for every even $k\ge 2$, either either $J$ has diwidth at most 
$16(9k)^2$ or $J$ embeds the $9k\times 9k$ diwall. Thus, to complete the proof of \ref{mainthm}, it suffices to prove the following two statements, the first of which is trivial.
\begin{thm}\label{diwidthtodartwidth}
The dartwidth of $G$ is at most the diwidth of $J$.
\end{thm}
\Proof
Let $(T,\tau)$ be a carving of $V(J)$ such that
for each
edge $e$ of $T$, the set of edges of $J$ joining the two parts of the
corresponding partition is a bond of $J$ with change number at most $k$. There is a natural bijection between the darts of $G$ and the vertices of $J$, so we can view $(T,\tau)$ as a carving of the set of darts of $G$; and as such, it shows that
the dart-width of $G$ is at most the diwidth of $J$, as required.~\bbox

\begin{thm}\label{keepwall2}
For all even $k\ge 2$, if $J$ embeds the $3k\times 9k$ diwall, then $G$ embeds the $k\times k$ diwall.
\end{thm}
\Proof  Suppose that $J$ embeds the $3k\times 9k$ diwall, and let $H$ be a subdigraph of $G$ that is a subdivision of 
this diwall. Let 
$P_1\LL P_{3k}$ and $Q_1\LL Q_{9k}$ be its horizontal and vertical paths, numbered in order.
Thus they form a $3k\times 9k$ diwall layout.
For each $v\in V(G)$, the cycle $C_v$ of $J$
bounds a region of $J$, and which is therefore a subset of one of the regions of the (unique, up to homeomorphism) didrawing of $H$. One of the 
regions of $H$ is the infinite region, but the boundaries of all its finite regions consist of subpaths of two consecutive of $P_1\LL P_{3k}$ and two consecutive of $Q_1\LL Q_{9k}$. Thus, for $v\in V(G)$, either $C_v$ intersects only 
some of $P_1,P_{3k}, Q_1,Q_{9k}$, or it intersects at most two consecutive members of $P_1\LL P_{3k}$ and two consecutive 
members of $Q_1\LL Q_{9k}$. Hence, the paths
$P_1,P_4,P_7\LL P_{3k-2}$ and $Q_1,Q_4,Q_7\LL Q_{9k-2}$ form a $k \times 3k$ diwall layout, and each $C_v$ intersects 
at most one of $P_1,P_4,P_7\LL P_{3k-2} $ and at most one of $Q_1,Q_4,Q_7\LL Q_{9k-2}$. 
The union of these paths is not a subdivision of a diwall, because (for instance) the last vertex of $P_1$
is not in any of $Q_1,Q_4,Q_7\LL Q_{9k-2}$; but by removing some vertices from the ends of these paths, we obtain
a  $k \times 3k$ diwall layout $R_1\LL R_{k}$ and $S_1\LL S_{3k}$ whose union is a $k \times 3k$ diwall $J$.

When we contract the edges of all the cycles $C_v$, we may no longer contain a $k \times 3k$ diwall layout; but we claim that we still
contain a $k \times 3k$ box system and hence a $k\times k$ diwall layout, by \ref{boxestodiwall}.
To see this, let us examine in more detail the effect on $R_1\LL R_{3k}$ and $S_1\LL S_{3k}$ of contracting the edges 
of some $C_v$. In any didrawing, if a region $r$ is incident with an edge $e$, we can speak of $r$ being to the left or 
right of $e$ in the natural sense (recall that we have assigned an orientation ``clockwise'' to $\Sigma$). In our case, 
if $C_v$ meets some $R_i$ or $S_j$, then it shares at least one edge with that path, since $G$ has maximum total degree 
three.

Now, $C_v$ bounds a region $r_v$ of $G$ that is a subset of a region $r'_v$ of the didrawing of $J$. Since
$C$ bounds $r_v$ in the clockwise direction from the construction, it follows that $r_v$ (and hence $r'_v$) is on the 
right of any edge of $J$ that belongs to $C_v$. Since the infinite region of $J$ is on the left of all edges of $J$
incident with it, $r'(v)$ is a finite region if $C_v$ shares any vertices at all with $R_1\LL R_{3k}$ or $S_1\LL S_{3k}$.

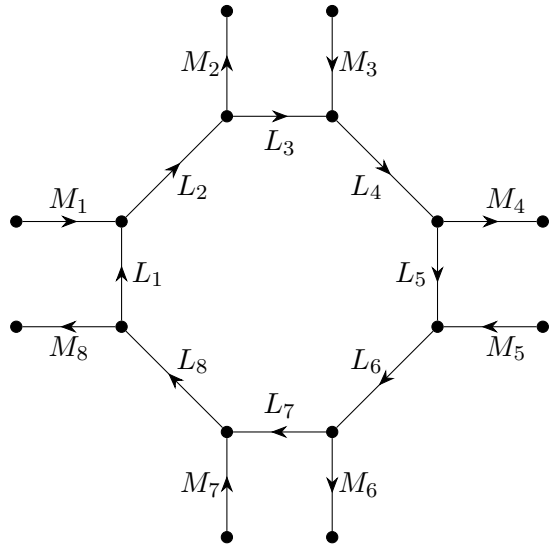
\begin{figure}[h!]
\centering

\begin{tikzpicture}[scale=.7,auto=left]

\tikzstyle{every node}=[inner sep=1.5pt, fill=black,circle,draw]

\node (1) at (-5,1) {};
\node (11) at (-3,1) {};
\node (8) at (-5,-1) {};
\node (18) at (-3,-1) {};
\node (4) at (5,1) {};
\node (14) at (3,1) {};
\node (5) at (5,-1) {};
\node (15) at (3,-1) {};
\node (2) at (-1,5) {};
\node (3) at (1,5) {};
\node (12) at (-1,3) {};
\node (13) at (1,3) {};
\node (7) at (-1,-5) {};
\node (17) at (-1,-3) {};
\node (6) at (1,-5) {};
\node (16) at (1,-3) {};
\foreach \from/\to in {
1/11, 12/2,3/13,14/4,5/15,16/6,7/17, 18/8, 
11/12,12/13,13/14,14/15,15/16,16/17,17/18,18/11}
\dedge{\from}{\to};

\tikzstyle{every node}=[]

\def\s{1.65}
\node at (-2.5,0) {$L_1$};
\node at (0,2.5) {$L_3$};
\node at (2.5,0) {$L_5$};
\node at (0,-2.5) {$L_7$};
\node at (-\s,\s) {$L_2$};
\node at (\s,\s) {$L_4$};
\node at (\s, -\s) {$L_6$};
\node at (-\s,-\s) {$L_8$};
\node at (-4,1.4) {$M_1$};
\node at (-1.5,4) {$M_2$};
\node at (1.5,4) {$M_3$};
\node at (4.3,1.4) {$M_4$};
\node at (4.3,-1.4) {$M_5$};
\node at (1.5,-4) {$M_6$};
\node at (-1.5, -4) {$M_7$};
\node at (-4,-1.4) {$M_8$};

\end{tikzpicture}
\caption{A piece of a subdivided diwall.} \label{fig:diwallpiece}
\end{figure}

There are four kinds of finite region in the didrawing $J$ (see Figure \ref{fig:diwall}); those bounded by a 
clockwise cycle (coming from a cycle of length
eight in the diwall before subdivision), those bounded by an anticlockwise cycle of ``unsubdivided length'' eight, and 
those bounded by cycles of unsubdivided length four (and there are two kinds of the last, depending whether the region is
to the right or left of the edges in horizontal paths incident with it). We can assume that $C_v$
shares an edge with at least one of $R_1\LL R_{3k}$ or $S_1\LL S_{3k}$, and so $r'(v)$ is not bounded by an 
anticlockwise cycle; and if $r'(v)$ is bounded by a cycle with unsubdivided length four, then it only intersects one of 
$R_1\LL R_{3k}$ or $S_1\LL S_{3k}$ and will not cause us a problem. The awkward case is when $r'(v)$
is bounded by a clockwise cycle of unsubdivided length eight (see Figure \ref{fig:diwallpiece}.)

Let $r'$ be some such region of $J$.
There are (in general) eight vertices of total degree three incident with $r'$, dividing its boundary into eight paths labelled as in the figure.

Let $r'$ have boundary formed by subpaths of
$R_i, R_{i+1}, S_j, S_{j+1}$.  Thus $i,j$ are both odd, since $r'$
is bounded by a clockwise cycle of unsubdivided length eight. In the figure, 
\begin{itemize}
\item $R_i$ contains $M_1, L_2,L_3,L_4, M_4$;
\item $R_{j+1}$ contains $M_3, L_4,L_5,L_6, M_6$;
\item $R_{i+1}$ contains $M_5, L_6, L_7, L_8, M_8$; and
\item $R_{j}$ contains $M_7,L_8, L_1,L_2,M_2$.
\end{itemize}
Let $X$ be the set of all $v\in V(\Gamma_{3,3k})$ with $r'(v) = r'$. 
For each $v\in X$, since $C_v$ does not intersect both $R_i, R_{i+1}$, and does not intersect both $S_j, S_{j+1}$, it follows that 
the set of  edges of $L_1\cupcup L_8$ in $C_v$ is included in one of $L_1\cup L_2\cup L_3$, $L_3\cup L_4\cup L_5$,
$L_5\cup L_6\cup L_7$, $L_7\cup L_8\cup L_1$. From planarity,  for $h = 1,3,5,7$  we can choose an edge
of $L_i$ that divides it into two 
vertex-disjoint
subpaths $L_h', L_h''$, such that  
$L_1''\cup L_2\cup L_3'$, $L_3''\cup L_4\cup L_5'$,
$L_5''\cup L_6\cup L_7'$, $L_7''\cup L_8\cup L_1'$ are directed paths, and 
for all $v\in X$, the set of  edges of $L_1\cupcup L_8$ in $C_v$ is included 
in one of these four paths.

For $1\le i\le k$ and $1\le j\le 3k$ with $i,j$ both odd, 
\begin{itemize}
\item let $A_{i,j}$ be the union of $V(L_1''\cup L_2\cup L_3')$ and all internal vertices of $M_1$;
\item let $A_{i,j+1}$ be the union of $V(L_3''\cup L_4\cup L_5')$ and all internal vertices of $M_5$; 
\item let $A_{i+1,j+1}$ be the union of $V(L_5''\cup L_6\cup L_7')$ and all internal vertices of $M_7$;
and 
\item let $A_{i+1,j}$ be the union of $V(L_7''\cup L_8\cup L_9')$ and all internal vertices of $M_9$.
\end{itemize}
This defines $A_{i,j}$ for all $i,j$ with $1\le i\le k$ and $1\le j\le 3k$ (with a few natural adjustments when $i = 1,k$ 
or $j=1,3k$ which we omit). Hence there is a $k\times 3k$ box system in $G$ where the vertex in the $i$th horizontal
path and $j$th vertical path of $\Gamma_{k,3k}$ is mapped to $A_{i,j}$. This box system has the property that for each
$v\in V(\Gamma_{k,3k})$, one of the boxes $A_{i,j}$ includes the set of all vertices in $C_v\cap J$, and so this remains 
a box system after contracting the edges of all the $C_v$. This proves \ref{keepwall2}, and hence proves \ref{mainthm}.~\bbox

\section{Grids}
We all agree what an undirected grid is, but there are several ways to define a directed grid. And while every undirected planar graph is a
minor of a big enough undirected grid, this is by no means true for planar digraphs and directed minor containment and 
directed grids. 

A $k\times k$  grid has $k$ ``horizontal'' paths and $k$ ``vertical'' paths, defined in the natural way, and
let us make them all directed paths,
but how should we direct them?

Three types that we care about are:
\begin{itemize}
\item All the horizontal paths are directed from left to right; we call this the $k\times k$ {\em acyclic} grid.
\item Assume $k$ is even, and the top $k/2$ horizontal paths are directed from left to right, and the others from right to left, and similarly the $k/2$ left-most vertical paths are directed from
bottom to top, and the others from top to bottom. (See Figure \ref{fig:semigrid}.) Let us call this the $k\times k$ {\em semi-grid}.
\item The $k\times k$ {\em alternating} grid, defined earlier.
\end{itemize}

For the acyclic grid, in a sense it makes no difference
how we direct the vertical paths, because for any acyclic grid $H$, if $k$ is large enough then every $k\times k$ acyclic grid
embeds $H$.
To see this, let $A_k$  be the $k\times k$ acyclic grid when all vertical paths are directed from top to bottom, and let
$B_k$ be the one when they
alternate in direction.
If $k'\ge 2k-1$
then clearly, every $k'\times k'$ acyclic grid embeds $A_k$; and
for every $k\times k$ acyclic grid $H$, if $k'\ge 2k-1$
then $B_{k'}$ embeds $H$.
But one can show that for
all $k$ there exists $\ell=O(k^2)$ such that $A_\ell$ embeds $B_k$.
Thus, up to embedding, these types of grid are all equivalent.
But we do not know any analogue of \ref{GM5} for acyclic grids. (That is an interesting question, perhaps: is there a version of \ref{GM5} for acylic digraphs?)

The semi-grid,
is equivalent (up to embedding, as before) to the grid when we
direct the horizontal paths as before (that is, the top $k/2$ horizontal paths are directed
from left to right, and the others from right to left), but have the vertical paths alternate in direction.
It is also equivalent in the same sense to what is called the $k\times 2k$ {\em cylindrical grid}, defined as follows.
Take $k$ disjoint cycles drawn in the plane, each inside the next and each rotating clockwise; and then take
$2k$ directed paths in the plane between the outside cycle and the inside one, alternating in direction and each intersecting each of the cycles in just one vertex. Again,
a large enough semi-grid embeds a large cylindrical grid, and vice versa.
For this one there is an analogue of \ref{GM5}, that we discuss next.

\section{Minor containment, and the Kawarabayashi-Kreutzer theorem}

Minor containment in digraphs is much more complicated than in undirected graphs. We agree what a minor
of an undirected graph means,
but for digraphs there are at least two definitions of minor containment of interest, ``butterfly'' minors and ``strong'' minors, and here we want a third, ``semi-strong'' minors.
\begin{itemize}
\item A {\em butterfly contraction} is the operation on a digraph of
contracting a non-loop edge $e=uv$ with the property that either $u$ has out-degree one or $v$ has in-degree one.
A digraph $H$ is a {\em butterfly minor} of a digraph $G$ if a digraph isomorphic to $H$ can be obtained
from a subdigraph of $G$ by a sequence of butterfly contractions.
\item A {\em strong contraction} is the operation on a digraph of contracting a non-null strongly-connected subdigraph to a single vertex.
A digraph $H$ is a {\em strong minor} of a digraph $G$ if a digraph isomorphic to $H$ can be obtained
from a subdigraph of $G$ by a sequence of strong contractions.
\item $H$ is a
 {\em semi-strong minor} of a digraph $G$ if some directed subdivision of $H$ is a strong minor
of $G$.
\end{itemize}

Not every didwawing is a strong minor of a diwall, or a butterfly minor of a diwall, but each is a semi-strong minor of a large eneough 
diwall (we omit the proof) and that is the reason we want semi-strong minors.
So a consequence of our main theorem is that:
\begin{thm}\label{strongthm}
For every didrawing
$H$ there exists $k$ such that every weakly 2-edge-connected loopless didrawing that does not contain $H$ as a semi-strong minor has dart-width at most $k$.
\end{thm}
It is easy to show the weak converse, that large enough diwalls have large dart-width, so our
result \ref{strongthm} is thus a sort of analogue of
the theorem of Robertson and Seymour~\cite{GM5} for planar digraphs.

There is at least one other way to extend that theorem to digraphs, due to Kawarabayashi and Kreutzer~\cite{kawa1, kawa2},
using butterfly minors and the semi-grid of Figure \ref{fig:semigrid}.
It needs ``directed tree-width'' which is a complicated parameter introduced in~\cite{johnson}. (We omit the definition.)
Kawarabayashi and Kreutzer~\cite{kawa1,kawa2} proved a conjecture
of~\cite{johnson}, that:
\begin{thm}\label{kawathm}
For every $k\ge 0$ there exists $\ell\ge 0$ such that every digraph with directed tree-width at least $\ell$ contains a
$k\times k$ semi-grid as a
butterfly minor. Conversely, for every $\ell\ge 0$ there exists $k\ge 0$ such that every digraph that contains a $k\times k$ semi-grid as a
butterfly minor has directed tree-width at least $\ell$.
\end{thm}

Consequently, we have the same conclusion if we exclude
a digraph $H$
that is a butterfly minor of a semi-grid grid. Such digraphs $H$ are the digraphs that
can be drawn in the plane without crossings,
such that every edge subtends a positive angle at the origin, and for every vertex $v$, the out-edges at $v$ and in-edges at $v$
fall into two disjoint intervals in the circular order around $v$ (we omit the proof). Let us call such digraphs {\em circular}. Thus, the
theorem of Kawarabayashi and Kreutzer says:
\begin{thm}\label{kawa}
For every circular digraph $H$, there exists $k$ such that every digraph that does not contain $H$ as a butterfly minor has directed tree-width at most $k$.
\end{thm}

Our result \ref{strongthm} says something similar to this,  but is different in three ways,  one gain, one loss and one tie:
\begin{itemize}
\item it involves semi-strong minors instead of butterfly minors (a tie);
\item the digraph $H$ we exclude can be any planar digraph $H$ (a gain);
\item it only tells us how to decompose {\em planar} digraphs that do not contain $H$ as  a semi-strong minor (a loss).
\end{itemize}

Finally, let us say a didrawing with interleaving at most two is {\em normal}. Thus, circular didrawings are normal.
The didrawings that are butterfly minors of diwalls are precisely the normal didrawings (again, we omit the proof)
so there is another consequence of our theorem:
\begin{thm}\label{butterflyversion}
For every normal didrawing $H$, there exists $k$ such that every didrawing not containing $H$
as a butterfly minor has dart-width at most $k$.
\end{thm}
\section*{Ackowledgement}
Part of this work was done during the 2025 Bertinoro Workshop on Algorithms and Graphs.
The authors are grateful to the organizers for the invitation, and to
the Bertinoro University Residential Center for its hospitality.  Thanks also to Alex Divoux for a careful reading.

\end{document}